\documentclass[twocolumn]{bmcart}

\usepackage{amsmath,amsthm,amssymb}

\usepackage{ascmac}
\usepackage{txfonts}
\usepackage[utf8]{inputenc} 

\usepackage{graphicx}
\usepackage{relsize}

\theoremstyle{definition}
\newtheorem{theorem}{Theorem}
\newtheorem{definition}{Definition}
\newtheorem{proposition}{Proposition}
\newtheorem{corollary}{Corollary}
\newtheorem{example}{Example}
\newtheorem{remark}{Remark}
\theoremstyle{remark}
\newtheorem{Proof}{Proof.}

\startlocaldefs
\endlocaldefs

\begin{document}

\begin{frontmatter}

\begin{fmbox}


\title{A Property of Random Walks on a Cycle Graph}


\author[
   addressref={aff1},                   
   email={y-ikeda(at)math.kyushu-u.ac.jp}   
]{\inits{YI}\fnm{Yuki} \snm{Ikeda}}

\author[
   addressref={aff2},
   email={fukai(at)math.kyushu-u.ac.jp}
]{\inits{YF}\fnm{Yasunari} \snm{Fukai}}

\author[
   addressref={aff3},
   email={ym(at)imi.kyushu-u.ac.jp}
]{\inits{YM}\fnm{Yoshihiro} \snm{Mizoguchi}}


\address[id=aff1]{
  \orgname{Graduate School of Mathematics, Kyushu University}, 
}
\address[id=aff2]{%
  \orgname{Faculty of Mathematics, Kyushu University},
}

\address[id=aff3]{%
  \orgname{Institute of Mathematics for Industry, Kyushu University},
}


\begin{artnotes}
\note[id=n1]{Equal contributor} 
\end{artnotes}



\begin{abstractbox}

\begin{abstract} 
We analyze the Hunter vs. Rabbit game on a graph, which is a model of communication in adhoc mobile networks.
Let $G$ be a cycle graph with $N$ nodes.
The hunter can move from a vertex to a vertex along an edge.
The rabbit can jump from any vertex to any vertex on the graph.
We formalize the game using the random walk framework.
The strategy of the rabbit is formalized using a one dimensional random walk over $\mathbb{Z}$. 
We classify strategies using the order $O(k^{-\beta-1})$ of their Fourier transformation.
We investigate lower bounds and upper bounds of the probability that the hunter catches the rabbit.
We found a constant lower bound if $\beta \in (0,1)$ which does not depend on the size $N$ of the graph.
We show the order is equivalent to $O(1/\log N)$ if $\beta=1$
and a lower bound is $1/N^{(\beta-1)/\beta}$ if $\beta \in (1,2]$.
These results help us to choose the parameter $\beta$ of 
a rabbit strategy according to the size $N$ of the given graph. 
We introduce a formalization of strategies using a random walk,
theoretical estimation of bounds of a probability that the hunter catches the rabbit, and also show computing simulation results.
\end{abstract}




\end{abstractbox}
\end{fmbox}

\end{frontmatter}



\section{Introduction}
We consider a game played by two players: the hunter and the rabbit.
This game is described using a graph $G(V,E)$ where $V$ is a set of vertices and $E$ is a set of edges.
Both players may use a randomized strategy.
The hunter can move from vertex to vertex along edges.
The rabbit can move to any vertex at once.
The hunter's purpose is to catch the rabbit in as few steps as possible.
On the other hand, the rabbit considers a strategy that maximizes the time until the hunter catch the rabbit.
If the hunter moves to a vertex that the rabbit is at, the game finishes and we say that the hunter catches the rabbit.
\par
The Hunter vs. Rabbit game model is used for analyzing  transmission procedures in mobile adhoc networks\cite{chatzigiannakis20011,chatzigiannakis20012}.
This model helps to send an electronic messages efficiently using mobile phones.
The expected value of time until the hunter catches the rabbit is equal to the expected time until the recipient receives the mail.
One of our goals is to improve these procedures.\par
We introduce some games resembling the Hunter vs. Rabbit game.
The first one is the Princess vs. Monster game.
In this game, the Monster tries to catch the Princess in area $D$.
The difference between the Hunter vs. Rabbit game is that the Monster catches the Princess if the distance between the two players is smaller than a chosen value.
Also the Monster moves at a constant speed whereas the Princess can move at any speed.
This game is played on a cycle graph as introduced by Isaacs\cite{isaacs1965}.
The Princess vs. Monster game has been investigated by Alpern \cite{alpern1974}, Zelikin \cite{zelikin1972}, and so on.
Gal analyzed the Princess-Monster game on a convex multidimensional domain \cite{gal1979}.\par
The next one is the Deterministic pursuit-evasion game.
In this game we consider 
a runaway hide dark spot, for example a tunnel.
Parsons innovated the search number of a graph\cite{parsons1976,parsons1978}.
The search number of a graph is the least number of people that are required to catch a runaway hiding dark spot moving at any speed.
LaPaugh \cite{lapaugh1993} showed that if the runaway is known not to be in edge $e$ at any point of time, then the runaway can not enter edge $e$ without being caught in the remainder of the game.
%
%
Meggido showed that the computation time of the search number of a graph is NP-hard\cite{megiddo1988}.
If an edge can be cleared without moving along it, but it suffices to 'look into' an edge from a vertex, 
then the minimum number of guards needed to catch the fugitive is called the node search number of graph \cite{kirousis1986}.
The pursuit evasion problem in the plane was introduced by Suzuki and Yamashita \cite{suzuki1992}.
They gave necessary and sufficient conditions for a simple polygon to be searchable by a single pursuer.
Later Guibas et al. \cite{guibas1999} presented a complete algorithm and showed that the problem of determining the minimal number of pursuers needed to clear a polygonal region with holes is NP-hard.
Park et al. \cite{park2001}
gave three necessary and sufficient conditions for a polygon to be searchable and showed that there is $O(n^2)$ time algorithm for constructing a search path
for an $n$-sided polygon.
Efrat et al. \cite{efrat2000} gave a polynomial time algorithm for the problem of clearing a simple polygon with a chain of $k$ pursuers when the first and 
last pursuer can only move on the boundary of the polygon.\par
A first study of the Hunter vs. Rabbit game can be found in \cite{aleliunas1979}. 
The presented hunter strategy is based on random walk on a graph and it is shown that the hunter catches an unrestricted rabbit within $O(nm^2)$ rounds, 
where $n$ and $m$ denote the number of nodes and edges, respectively.
Adler et al. showed that if the hunter chooses a good strategy, the upper bound of the expected time that the hunter catches the rabbit is $O(n\log (diam(G)))$, where $diam(G)$ is a diameter of a graph $G$, and 
if the rabbit chooses a good strategy, the lower bound of the expected time that the hunter catches the rabbit is $\Omega(n\log (diam(G)))$ \cite{adler2003}.
Babichenko et al. showed Adler's strategies yield a Kakeya set consisting of $4n$ triangles with minimal area \cite{babichenko2012}.\par
In this paper, we propose three assumptions for the strategy of the rabbit.
We have the general lower bound formula for the probability that the hunter catches the rabbit.
The strategy of the rabbit is formalized using a one dimensional random walk over $\mathbb{Z}$. 
We classify strategies using the order $O(k^{-\beta-1})$ of their Fourier transform.
If $\beta = 1$, the lower bound of a probability that the hunter catches the rabbit is $((c_{*}\pi)^{-1}\log N + c_2)^{-1}$ where 
$c_{2}$ and $c_{*}$ are constants defined by the given strategy.
If $\beta \in (1,2]$, the lower bound of the probability that the hunter catches the rabbit is $c_{4}N^{-(\beta-1)/\beta}$ where $c_4 > 0$ is are constant defined by the given strategy.\par
We show experimental results for three examples of the rabbit strategy.
\begin{enumerate}
\item $\displaystyle
P\left\{X_{1}=k\right\} =\left\{
\begin{array}{ll}
\displaystyle\frac{1}{2a(\vert k\vert+1)(\vert k\vert+2)} &\quad (k\in \mathbb{Z}
 \setminus\left\{0\right\})\\
\displaystyle 1-\frac{1}{2a} &\quad (k=0)
\end{array}
\right.
$
\item $\displaystyle
P\left\{X_{1}=k\right\}=\left\{
\begin{array}{ll}
\displaystyle\frac{1}{2a\vert k\vert^{\beta +1}} &\quad (k\in\mathbb{Z}\setminus\left\{0\right\})\\
\displaystyle 1-\frac{1}{a}\sum_{k=1}^{\infty}\frac{1}{k^{\beta +1}} &\quad (k=0)
\end{array}
\right.
$
\item $\displaystyle
P\left\{X_{1}=k\right\}=\left\{
\begin{array}{ll}
\displaystyle\frac{1}{3} &\quad (k\in\left\{-1,0,1\right\})\\
\displaystyle 0 &\quad (k\not\in\left\{-1,0,1\right\}).
\end{array}
\right.
$
\end{enumerate}
We can confirm our bounds formula, and the asymptotic behavior of those bounds by the results of simulations.

\section{Statements of Results}
We consider the Hunter vs Rabbit game on a cycle graph.
To explain the Hunter vs Rabbit game, we  introduce some notation. 
Let $ X_{1}, X_{2}, \ldots $ be independent, identically distributed random variables defined on a probability space $(\Omega, {\cal F}, P)$ taking values in the integer lattice ${\mathbb{Z}}$. 
A one-dimensional random walk $ \{ S_{n} \}_{n=1}^{ \infty }$ is defined by 
$$ S_{n}= \sum_{j=1}^{n} X_{j}.$$
Let $ Y_{1}, Y_{2}, \ldots $ be independent, identically 
distributed random variables defined on a probability space $(\Omega_{\mathcal{H}}, {\cal F}_{\mathcal{H}}, P_{\mathcal{H}})$ taking values in the integer lattice ${\mathbb{Z}}$ with 
$$ P_{\mathcal{H}} \{ \vert Y_{1} \vert \leq 1 \} =1. $$
Let $ N \in {\mathbb{N}}$ be fixed. 
We denote by $X_{0}^{(N)} $ a random variable defined on a probability space $(\Omega_{N}, {\cal F}_{N}, \mu_{N})$ taking values in $ V_{N}:= \{ 0,1,2, \ldots , N-1 \}$ with 
$$ \mu_{N} \{ X_{0}^{(N)}=l \} = \frac{1}{N} \quad ( l \in V_{N}) .$$

For $ b \in {\mathbb{Z}}$, we denote by $(b \mod N)$ the remainder of $b$ divided by $N$. 

A rabbit's strategy $ \{ \mathcal{R}_{n}^{(N)} \}_{n=0}^{ \infty } $ is defined by 
$$  \mathcal{R}_{0}^{(N)}  = X_{0}^{(N)} \ \mbox{ and } \ 
\mathcal{R}_{n}^{(N)} =( X_{0}^{(N)} + S_{n} \mod N) .$$ 
$ \mathcal{R}_{n}^{(N)} $ indicates the position of the rabbit at time $n$ on $V_{N}$. 
Hunter's strategy $ \{ \mathcal{H}_{n}^{(N)}  \}_{n=0}^{ \infty }  $ is defined by 
$$  \mathcal{H}_{0}^{(N)} =0 \ \mbox{ and }  \mathcal{H}_{n}^{(N)}  =\left( \sum_{j=1}^{n} Y_{j} \mod N\right). $$
$ \mathcal{H}_{n}^{(N)} $ indicates the position of the hunter at time $n$ on $V_{N}$. 
Put 
$$  \mathbb{P}_{\mathcal{R}}^{(N)} = \mu_{N} \times P \ \mbox{ and } 
\ \tilde{ \mathbb{P} }^{(N)}= P_{\mathcal{H}} \times  \mathbb{P}_{\mathcal{R}}^{(N)} . $$
The hunter catches the rabbit when
the hunter and the rabbit are both located on the same place.

We will discuss the probability that the hunter catches the rabbit by time $N$ on $V_{N}$, that is,  
$$ \tilde{ \mathbb{P} }^{(N)} \left ( \bigcup_{n=1}^{N} \{ 
 \mathcal{H}_{n}^{(N)} =
 \mathcal{R}_{n}^{(N)}   \} \right) .$$
We investigate the asymptotic estimate of this probability 
as $ N \rightarrow \infty $.

\begin{definition}
We define conditions (A1), (A2) and (A3) as follows.
\begin{itemize}
\item[$(A1)$] The random walk $\{S_{n} \}_{n=1}^{\infty}$ is strongly aperiodic, i.e. 
for each $y \in \mathbb{Z}$, the smallest subgroup containing the set 
\begin{eqnarray*}
 \left\{y+k \in  {\mathbb{Z}} \ \vert \ P\left\{X_1 = k \right\} > 0\right\}
\end{eqnarray*}
is $\mathbb{Z}$.
\item[$(A2)$] $P\left\{X_{1}= k \right\} = P\left\{X_1=- k  \right\} 
\quad (k \in 
{\mathbb{Z}})$.
\item[$(A3)$] There exist $\beta \in (0,2]$, $c_{*}>0$ and  $\varepsilon >0$ 
such that
\begin{eqnarray*}
\phi(\theta) := \sum_{k \in \mathbb{Z}}e^{i\theta k}P\left\{X_1 = k \right\} 
= 1 - c_{*} \vert \theta \vert ^{\beta} + O(\vert \theta \vert ^{\beta + \varepsilon}).
\end{eqnarray*}
\end{itemize}
We denote the $\beta$ in $(A3)$ as $\beta_{\mathcal{R}}$.
\end{definition}

\begin{theorem}\label{main_theorem}
Assume that $X_{1}$ satisfies $(A1)-(A3)$.
\begin{enumerate}
\item [(I)]
If $\beta_{\mathcal{R}} \in (0,1)$, then there exists a constant $c_{1}>0$ 
such that for $ N \in \mathbb{N} \setminus \{ 1 \} $ and $  y_{1} ,y_{2}, \ldots  ,
y_{N} \in {\mathbb{Z}}$ with $ \vert y_{n}-y_{n+1} \vert \leq 1 \ (
n=1,2, \ldots , N-1)$, 
\begin{equation} 
c_{1} \leq  \mathbb{P}_{\mathcal{R}}^{(N)} \left( \bigcup_{n=1}^{N}
\left\{ \mathcal{R}_{n}^{(N)} = ( y_{n} \mod N) \right\} \right) . 
\label{eq:r1}
\end{equation}

\item[(I\hspace{-.1em}I)]
If $\beta_{\mathcal{R}}=1$, then there exist constants 
$c_{2}>0$ and $c_{3}>0$ 
such that for $ N \in \mathbb{N} \setminus \{ 1 \} $ and $  y_{1} , y_{2}, \ldots  ,
y_{N} \in {\mathbb{Z}}$ with $ \vert y_{n}-y_{n+1} \vert \leq 1 \ (
n=1,2, \ldots , N-1)$,
\begin{eqnarray}
\frac{1}{ \frac{1}{ c_{*} \pi} \log N +c_{2} } &\leq & \mathbb{P}_{\mathcal{R}}^{(N)} \left( \bigcup_{n=1}^{N}\left\{ \mathcal{R}_{n}^{(N)} = ( y_{n} \mod N) \right\} \right)\nonumber\\
&\leq & \frac{ c_{3}}{ \log N} .\label{eq:r2}
\end{eqnarray}

\item[(I\hspace{-.1em}I\hspace{-.1em}I)]
If $\beta_{\mathcal{R}} \in (1,2]$, then there exists a constant $c_{4}>0$  
such that for $ N \in \mathbb{N} \setminus \{ 1 \} $ and $ y_{1} , y_{2}, \ldots  ,
y_{N} \in {\mathbb{Z}}$ with $ \vert y_{n}-y_{n+1} \vert \leq 1 \ (
n=1,2, \ldots , N-1)$,
\begin{eqnarray} 
\frac{ c_{4}}{ N^{( \beta-1)/ \beta}} \leq \mathbb{P}_{\mathcal{R}}^{(N)} \left( \bigcup_{n=1}^{N}\left\{ \mathcal{R}_{n}^{(N)} = ( y_{n} \mod N) \right\} \right).
\label{eq:r3}
\end{eqnarray}
\end{enumerate}
\end{theorem}

The following bounds are obtained as a corollary of Theorem 1.

\begin{corollary}
Assume $(A1)-(A3)$.

If  $\beta_{\mathcal{R}} \in (0,1)$, then there exists a constant $c_{1}>0$ 
such that for $ N \in \mathbb{N} \setminus \{ 1 \} $, 
$$ c_{1} \leq \tilde{ \mathbb{P} }^{(N)} \left ( \bigcup_{n=1}^{N} \{ 
 \mathcal{H}_{n}^{(N)} =
 \mathcal{R}_{n}^{(N)}   \} \right). $$

If $\beta_{\mathcal{R}}=1$, then there exist constants 
$c_{2}>0$ and $c_{3}>0$ 
such that for $ N \in \mathbb{N} \setminus \{ 1 \} $, 
\begin{eqnarray}
\frac{1}{ \frac{1}{c_{*} \pi } \log N + c_{2}} &\leq & \tilde{ \mathbb{P} }^{(N)} \left ( \bigcup_{n=1}^{N} \left\{ \mathcal{H}_{n}^{(N)} = \mathcal{R}_{n}^{(N)}\right\} \right)\nonumber\\
&\leq & \frac{ c_{3}}{ \log N}. 
\label{eq:4}
\end{eqnarray}

If $\beta_{\mathcal{R}} \in (1,2]$, then there exists a constant $c_{4}>0$  
such that for $ N \in \mathbb{N} \setminus \{ 1 \} $, 
$$ \frac{ c_{4}}{ N^{( \beta-1)/ \beta}} \leq \tilde{ \mathbb{P} }^{(N)} \left ( \bigcup_{n=1}^{N} \{ 
 \mathcal{H}_{n}^{(N)} =
 \mathcal{R}_{n}^{(N)}   \} \right). $$
\end{corollary}

\begin{remark}\label{remark1}
Adler, R$\ddot{\mbox{a}}$cke, Sivadasan, Sohler and V$\ddot{\mbox{o}}$cking considered 
$ \tilde{ \mathbb{P} }^{(N)} ( \cup_{n=1}^{N} \{ 
 \mathcal{H}_{n}^{(N)} =
 \mathcal{R}_{n}^{(N)}   \})$
 in the case of
\begin{eqnarray*}
P\left\{X_{1}=k\right\} =\left\{
\begin{array}{ll}
\displaystyle\frac{1}{2(\vert k\vert+1)(\vert k\vert+2)} &\quad ( k \in \mathbb{Z} \setminus \left\{0\right\} ) \\
\displaystyle \frac{1}{2} &\quad (k=0).
\end{array}
\right.
\end{eqnarray*}
In this case, $ X_{1} $ satisfies $(A1)$, $(A2)$ and 
$$ \phi ( \theta )= 1 - \frac{ \pi }{2} \vert \theta \vert + 
O( \vert \theta \vert^{3/2}) $$
($(A3)$ with $ \beta =1 $), and we have (\ref{eq:4}) in Corollary 1 
which coincides with the result of Lemma 3 in  \cite{adler2003}. 
\end{remark}

\begin{remark}\label{remark2}
For $ \beta \in (0,2)$, let 
\begin{eqnarray*}
P \left\{X_{1}=k\right\}=\left\{
\begin{array}{ll}
\displaystyle\frac{1}{2a\vert k\vert^{\beta +1}} &\quad ( k\in {\mathbb{Z}} \setminus \left\{0\right\} ) \\
\displaystyle 1-\frac{1}{a}\sum_{k=1}^{\infty}\frac{1}{k^{\beta +1}} & \quad (k=0)
\end{array}
\right.
\end{eqnarray*}
with a constant $a$ satisfying $a > \sum_{k=1}^{ \infty} (1/ k^{ \beta +1}). $
Then $\phi(\theta)$ in $(A3)$ is 
\begin{equation}
  \phi ( \theta )= 1 - \frac{ \pi }{2a}
\frac{ \vert \theta \vert^{ \beta }}{ \Gamma ( \beta +1) \sin ( \beta \pi /2)} + 
O( \vert \theta \vert^{\beta +(2- \beta )/2}), \label{eq:AAAA}
\end{equation}
where $ \Gamma $ is the gamma function (see Appendix (B)). 
$X_{1}$ satisfies  $(A1)$, $(A2)$ and (\ref{eq:AAAA}).

If $X_{1}$ takes three values $ -1, 0,1$ with equal probability, then 
$X_{1}$ satisfies  $(A1)$, $(A2)$ and 
$$ \phi ( \theta )= 1- \frac{1}{3} \vert \theta \vert^{2} + O( \vert \theta \vert^{4}) $$
($(A3)$ with $ \beta =2$). 
\end{remark}

The inequality (\ref{eq:r3}) seems to be sharp, because the powers of upper and lower bound appearing in  (\ref{eq:r3}) cannot be improved. Indeed, we have the following estimates.

\begin{proposition}\label{prop;add1}
Let $\mathcal{H}_{i}^{(N)} = 0$ for any $i$ and assume $(A1)-(A3)$. If $\beta_{\mathcal{R}} \in (1,2]$, then there exist constants $c_{5},c_{6}>0$  
such that for $ N \in \mathbb{N}$, 
\begin{equation}
\frac{ c_{5}}{ N^{( \beta-1)/ \beta}} \leq 
\mathbb{P}_{\mathcal{R}}^{(N)} \left( \bigcup_{n=1}^{N}
\{ \mathcal{R}_{n}^{(N)} = 0  \} \right) 
\leq \frac{ c_{6}}{ N^{( \beta-1)/ \beta}}. \label{eq:00}
\end{equation}
\end{proposition}

\begin{proposition}\label{prop;add2}
Let $\mathcal{H}_{i}^{(N)} = i$ for any $i$.
If $X_{1}$ takes three values $ -1, 0,1$ with equal probability, 
then there exists a constant $c_{7}>0$  
such that for $ N \in \mathbb{N}$, 
\begin{equation}
c_{7}  \leq 
\mathbb{P}_{\mathcal{R}}^{(N)} \left( \bigcup_{n=1}^{N}
\{ \mathcal{R}_{n}^{(N)} = (n \mod N)  \} \right) .  \label{eq:000}
\end{equation}
\end{proposition}

The proofs of Proposition \ref{prop;add1} and Proposition \ref{prop;add2} are given in Appendix (A).

\begin{remark}\label{remark4}
Assume $(A1)$ and $(A2)$. If  there exist $c_{*}>0$ and  $\varepsilon >0$
such that
$$
\phi(\theta) 
= 1 - c_{*} \vert \theta \vert + O(\vert \theta \vert ^{1 + \varepsilon})
$$
($(A3)$ with $ \beta =1$). Then 
\begin{equation}
\lim_{ N \rightarrow \infty } \left( \frac{1}{ c_{*} \pi } 
\log N \right) \mathbb{P}_{\mathcal{R}}^{(N)} \left( \bigcup_{n=1}^{N}
\{ \mathcal{R}_{n}^{(N)} = 0  \} \right) =1. \label{eq:LIM}
\end{equation}

The proof of (\ref{eq:LIM}) is given in Appendix (C).
\end{remark}

\section{Computer simulation}
In this section, we show some experimental results about the Hunter vs Rabbit game on a cycle graph.
We compute $P\left\{S_n \mod N = k\right\}$ by using the gamma function and the class {\tt discrete\_distribution} in {\tt C++}.
We can show the probability the rabbit is caught and the expected value of the time until the rabbit is caught using this application.\par
In this section, we consider a lower bound $L(N,a)$ of the probability that the hunter catches the rabbit.
According to the Proposition~\ref{bound;thm2} and Proposition~\ref{proof;prop4}, we define $L(N,a)$ as follows: 
\begin{eqnarray*}
L(N) = \frac{1}{1 + A_{N} + B_{N} + \frac{1}{1-\rho_{*}}}
\end{eqnarray*}
where
\begin{eqnarray*}
A_{N} = \left\{
\begin{array}{ll}
\frac{2^{2+\varepsilon - \beta}\pi^{\varepsilon - \beta}C_{*}}{c_{*}^{2}} &\quad (\beta \in (0,1]),\\
2N^{(\beta -1)/\beta} &\quad (\beta \in (1,2))
\end{array}
\right.
\end{eqnarray*}
and
\begin{eqnarray*}
B_{N} = \left\{
\begin{array}{ll}
\frac{2^{1-\beta}}{\pi^{\beta}c_{*}(1-\beta)} &\quad (\beta \in (0,1)),\\
\frac{1}{\pi c_{*}}\log N + \frac{1}{\pi c_{*}} &\quad (\beta = 1),\\
\frac{2^{2-\beta}}{c_{*} \pi}\left(1+\frac{1}{\beta -1}\right)N^{(\beta -1)/\beta} &\quad (\beta \in (1,2)).
\end{array}
\right.
\end{eqnarray*}
We note $\beta$ and $c_*$ are defined by a given $P\{X_t=k\}$ in an example.
We choose appropriate constants $\varepsilon$, $\rho_*$ and  $C_*$  for each examples. 

\begin{example}\label{example1}
We consider the generalization of the case of \cite{adler2003}.
Let
\begin{eqnarray*}
P\left\{X_{t}=k\right\} =\left\{
\begin{array}{ll}
\displaystyle\frac{1}{2a(\vert k\vert+1)(\vert k\vert+2)} &\quad (k\in \mathbb{Z}
 \setminus\left\{0\right\})\\
\displaystyle 1-\frac{1}{2a} &\quad (k=0)
\end{array}
\right.
\end{eqnarray*}
where $a \geq \frac{1}{2}$.
We note $\beta=1$, $c_* = \pi$ and $\varepsilon=1/2$ in Remark~\ref{remark1}.
If $a=1$, then this is the case in \cite{adler2003}.
We can define $C_{*}$ and $\rho_{*}$ for this case.
So we have
\begin{eqnarray}
\frac{1}{\sum_{i=0}^{N-1}p_{i}^{(N)}} \geq L(N,1) = \frac{1}{\frac{2}{\pi^{2}}\log N + 6.50503}.
\label{comsi;ex1;eq0001}
\end{eqnarray}
The proof of (\ref{comsi;ex1;eq0001}) is given in Appendix (D).\par

Figure \ref{heavy_tailed} shows an experimental result of the probabilities for all initial positions of the rabbit with $N=100$ and $a = 1$.
The horizontal axis is the initial position of the rabbit,
and the vertical axis shows the probability the rabbit is caught.
The red line in the figure is a probability that the hunter catches the rabbit.
The blue line is the average of probabilities that the hunter catches the rabbit.
The green line is $L(N,a)$.
In this case, the hunter does not move from the initial position $0$.
As you can see, the average of the probability that the hunter catches the rabbit is bounded below by $L(N,a)$.
\par
In this case, the average of the probability that the hunter catches the rabbit each initial position of the rabbit nearly equals $0.4258$, so we have
\begin{eqnarray*}
\frac{1}{L(100,1)} \fallingdotseq 7.43823,
\end{eqnarray*}
and
\begin{eqnarray*}
\frac{1}{L(100,1)} \mathbb{P}_{\mathcal{R}}^{(N)} \left( \bigcup_{n=1}^{N}\{ \mathcal{R}_{n}^{(N)} = 0  \} \right) \fallingdotseq 3.1672.
\end{eqnarray*}

Table \ref{table;heavy_tailed} is the experimental results of Example \ref{example1} with $a=1$ and $N = 100, 500$ and  $1000$.
This table shows the asymptotic behavior of (\ref{eq:LIM}).

\begin{table}[htbp]
\label{table;heavy_tailed}
\caption{This table is experimental results of Example \ref{example1} with $a=1$ and $N = 100, 500$ and $1000$.
$A$ is the average of the probability that the hunter catches the rabbit.}
\begin{tabular}{c|ccc}
\hline
$N$    & $1/L(N,a)$  & $A$ & $A/L(N,a)$\\ \hline
$100$  & 7.43823 & 0.4528 & 3.1672\\
$500$  & 7.76437 & 0.39048 & 3.03183\\
$1000$ & 7.90483 & 0.37555 & 2.96866\\ \hline
\end{tabular}
\end{table}



\end{example}

\begin{example}\label{example2}
We consider the case of $\beta \in (0,2)$. We put
\begin{eqnarray*}
P\left\{X_{t}=k\right\}=\left\{
\begin{array}{ll}
\displaystyle\frac{1}{2a\vert k\vert^{\beta +1}} &\quad (k\in\mathbb{Z}\setminus\left\{0\right\})\\
\displaystyle 1-\frac{1}{a}\sum_{k=1}^{\infty}\frac{1}{k^{\beta +1}} &\quad (k=0)
\end{array}
\right.
\end{eqnarray*}
where $a > \sum_{k=1}^{\infty}\frac{1}{k^{\beta +1}}$.
By Remark \ref{remark2}, $c_{*} = \frac{\pi}{2a \Gamma ( \beta +1) \sin ( \beta \pi /2)}$ and $\varepsilon = \frac{2- \beta}{2}$.
Then, the lower bound of the probability that the hunter catches the rabbit $L(N,a)$ is
\begin{eqnarray*}
&&L(N,a)\\
&&= \left\{
\begin{array}{l}
\frac{1}{1+2^{1-\beta}(1-\beta)^{-1}\pi^{-\beta}c_{*}^{-1} + 2^{4-3\beta /2}\pi^{1-3\beta /2}c_{*}^{-1}C_{*} + (1-\rho_{*})^{-1}}\\
 \hspace{5cm} (\beta \in (0,1))\\
\frac{1}{1+(\pi c_{*})^{-1}(1+\log N)+2^{7/2}\pi^{-1/2}c_{*}^{-1}C_{*} + (1-\rho_{*})^{-1} }\\
 \hspace{5cm} (\beta =1)\\
\frac{1}{1+2N^{(\beta-1)/\beta}+2^{2-\beta}c_{*}^{-1}\pi^{-\beta}\left(1+(\beta-1)^{-1}
\right)N^{(\beta-1)/\beta}+(1-\rho_{*})^{-1}}\\
 \hspace{5cm} (\beta \in (1,2))
\end{array}
\right.
\end{eqnarray*}
where $\rho_*$ and  $C_*$ are appropriate constants for each examples. 
When $a=2.5$ and $\beta=1$, we set $C_* \fallingdotseq 0.177245$ and $\rho_{*} \fallingdotseq 0.694811$.
So we have
\begin{eqnarray*}
L(N,2.5) = \frac{1}{\frac{5}{\pi^{2}}\log N + 4.65936}.
\end{eqnarray*}
\par
Figure \ref{inverse} is an experimental result with $\beta = 1$, $N = 100$ and $a = 2.5$.
In this case, the average of the probability that the hunter catches the rabbit
 nearly equals $0.318$,  so we have
\begin{eqnarray*}
\frac{1}{L(100,2.5)} \fallingdotseq 6.99237,
\end{eqnarray*}
and
\begin{eqnarray*}
\frac{1}{L(100,2)} \mathbb{P}_{\mathcal{R}}^{(N)} \left( \bigcup_{n=1}^{N}\{ \mathcal{R}_{n}^{(N)} = 0  \} \right) \fallingdotseq 2.22357.
\end{eqnarray*}

Table \ref{table;inverse} is the experimental results of Example \ref{example2} with $\beta=1$, $a=2.5$ and $N = 100, 500$ and $1000$.
This table shows that the value of $A/L(N,a)(>1)$ is decreasing.
\begin{table}[htbp]
\label{table;inverse}
\caption{This table is experimental results of Example \ref{example2} with $\beta=1$, $a=2.5$ and $N = 100, 500$ and $1000$.
$A$ is the average of the probability that the hunter catches the rabbit.}
\begin{tabular}{c|ccc}
\hline
$N$    & $1/L(N,a)$  & $A$ & $A/L(N,a)$\\ \hline
$100$  & 6.99237 & 0.318   & 2.22357\\
$500$  & 7.80772 & 0.25924 & 2.02407\\
$1000$ & 8.15887 & 0.24015 & 1.95935\\ \hline
\end{tabular}
\end{table}


\begin{figure*}[htbp]
\begin{center}
\includegraphics[width=6cm,bb=5 5 330 272]{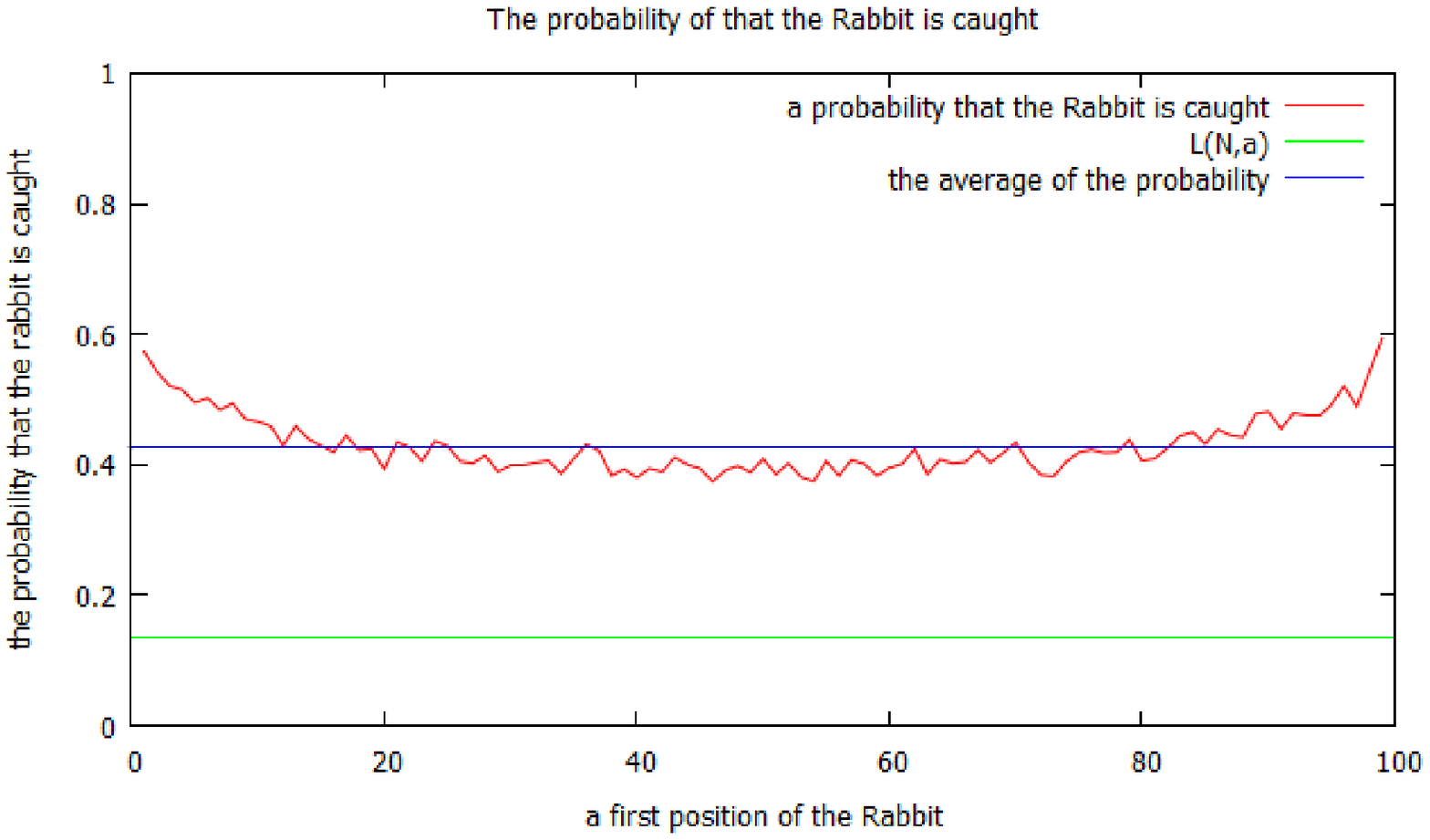}
\end{center}
\caption{
This is an experimental result of Example \ref{example1}.
In this case, $a=1$.
The hunter does not move from an initial position $0$.
}
\label{heavy_tailed}
\end{figure*}

\begin{figure*}[htbp]
\begin{center}
\includegraphics[width=6cm,bb=5 5 330 272]{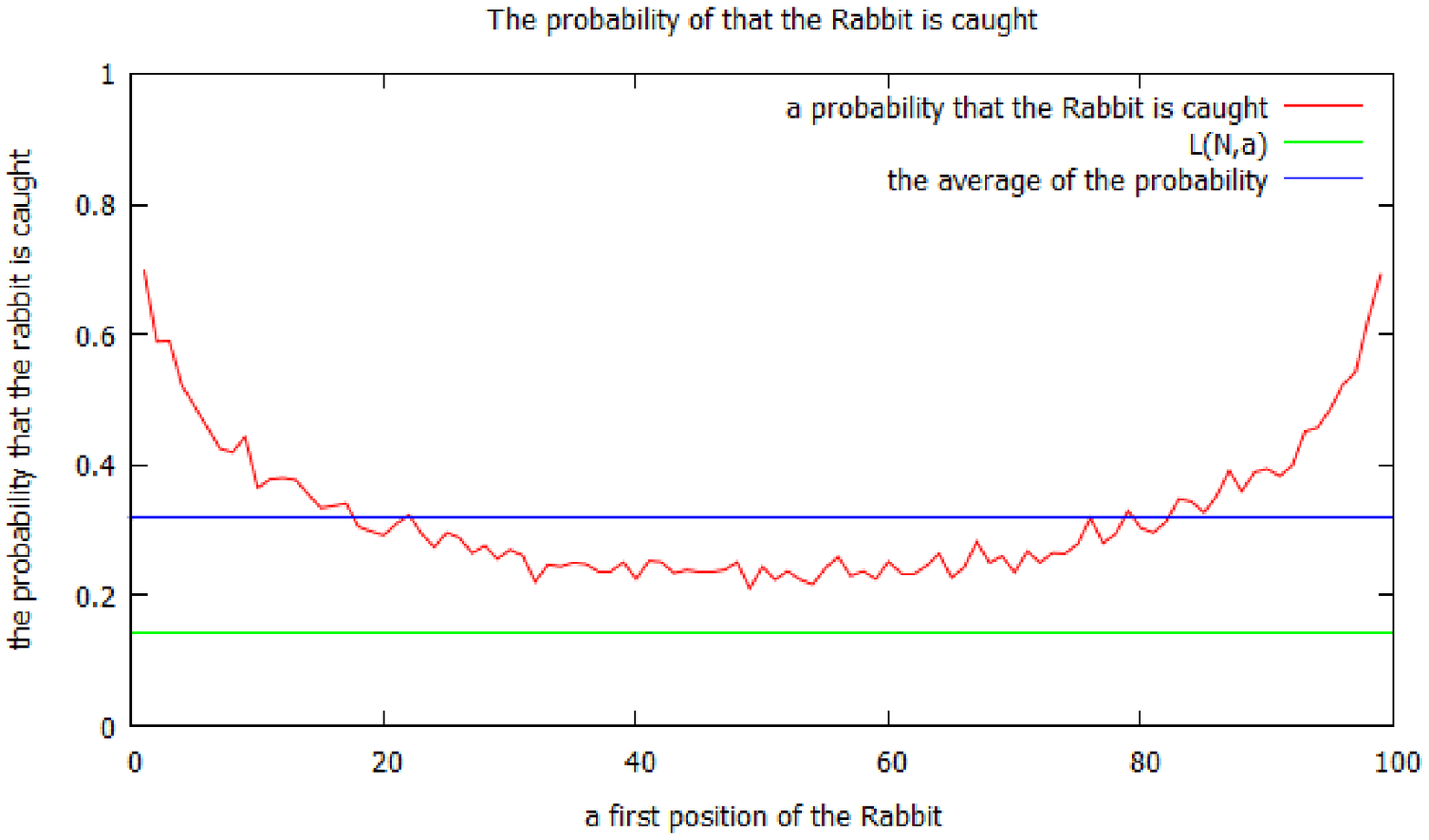}
\end{center}
\caption{
This is an experimental result of Example \ref{example2}.
In this case, $a=2.5$.
The hunter does not move from an initial position $0$.
}
\label{inverse}
\end{figure*}

\begin{figure*}[htbp]
\begin{center}
\includegraphics[width=6cm,bb=5 5 330 272]{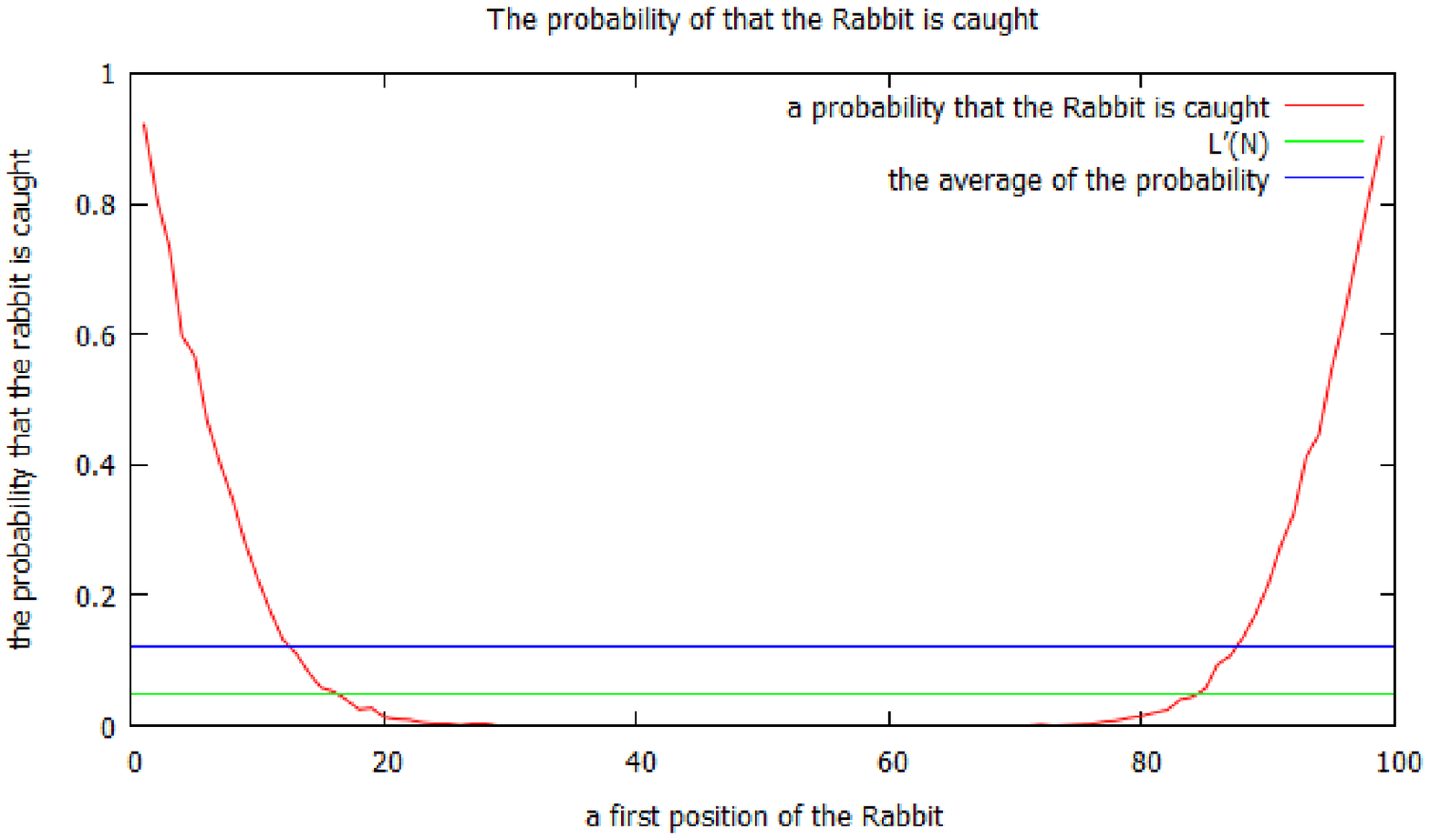}
\end{center}
\caption{
This is an experimental result of Example \ref{example3}.
The hunter does not move from an initial position $0$.
}
\label{slow_strategy}
\end{figure*}

\end{example}

\begin{example}\label{example3}
We put
\begin{eqnarray*}
P\left\{X_{t}=k\right\}=\left\{
\begin{array}{ll}
\displaystyle\frac{1}{3} &\quad (k\in\left\{-1,0,1\right\})\\
\displaystyle 0 &\quad (k\not\in\left\{-1,0,1\right\}).
\end{array}
\right.
\end{eqnarray*}
By Remark \ref{remark2}, $\beta = 2$, $c_{*} = \frac{1}{3}$ and $\varepsilon = 2$.
In this case, the lower bound of the probability the hunter catches the rabbit $L'(N)$ is
\begin{eqnarray*}
L'(N) = \frac{1}{\left(1+\frac{6}{\pi^2}\right)N^{1/2} + 4.26301}.
\end{eqnarray*}
(We can prove this using in the same way in Appendix (D).)
Figure \ref{slow_strategy} is an experimental result of Example \ref{example3}.
The green line in Figure \ref{slow_strategy} is $L'(N)$.
\end{example}
We could have a concrete lower bound of the average of a probability that the hunter catches the rabbit for those examples.

\section{Upper bounds and Lower bounds}

In this section, we give a relation between  
$$
\mathbb{P}_{\mathcal{R}}^{(N)} \left( \bigcup_{n=1}^{N}\left\{ \mathcal{R}_{n}^{(N)} = ( y_{n} \mod N) \right\} \right)
$$
and one-dimensional random walk $ \{ S_{n} \}_{n=1}^{ \infty}$.

\begin{proposition}\label{bound;thm2}
For $ N \in {\mathbb{N}} \setminus \{ 1 \}$ and $y_{1},y_{2}, \ldots , y_{N} \in {\mathbb{Z}} $ 
with $ \vert y_{n} -y_{n+1} \vert \leq 1 \ (n=1,2, \ldots , N-1)$, 
\begin{eqnarray}
\frac{1}{ \sum_{i=0}^{N-1} p_{i}^{(N)} } &\leq & \mathbb{P}_{\mathcal{R}}^{(N)} \left( \bigcup_{n=1}^{N}\left\{ \mathcal{R}_{n}^{(N)} = ( y_{n} \mod N) \right\} \right)\nonumber\\
&\leq &  \frac{2}{ \sum_{i=0}^{N-1} q_{i}^{(N)} },\label{eq:PNR}
\end{eqnarray}
where 
 $$ [y]_{N}= \{ y+kN \ | \ k \in {\mathbb{Z}} \} , $$
\begin{eqnarray*}
p_{i}^{(N)} =
\left\{
\begin{array}{ll}
1 &\quad (i=0) \\
\displaystyle \max_{\vert y \vert \leq i, \ y \in {\mathbb{Z}}} P\left\{S_{i} \in [y]_{N} \right\} &\quad (i \in \mathbb{N}) \\ 
\end{array}
\right. 
\end{eqnarray*}
and
\begin{eqnarray*}
q_{i}^{(N)} =
\left\{
\begin{array}{ll}
1 &\quad (i=0) \\
\displaystyle \min_{\vert y \vert \leq i, \ y \in {\mathbb{Z}}} P\left\{S_{i} \in [y]_{N} \right\} &\quad (i \in \mathbb{N}).
\end{array}
\right. 
\end{eqnarray*}
\end{proposition} 

\begin{Proof}

We note that
\begin{eqnarray*}
&&\bigcup_{n=1}^{N}\left\{ \mathcal{R}_{n}^{(N)} = ( y_{n} \mod N) \right\}\\
&&= \bigcup_{l=0}^{ N-1} \bigcup_{n=1}^{N} \left\{ X_{0}^{(N)}=l, \ l+S_{n} \in [y_{n}]_{N} \right\}\\
&&= \bigcup_{l=0}^{ N-1} \bigcup_{n=1}^{N}  \left\{ 
\begin{array}{ll}
X_{0}^{(N)}=l, &    
l+S_{n} \in [y_{n}]_{N} ,\\
 l+S_{i} \notin [y_{i}]_{N},  &
1 \leq i \leq n-1
\end{array} 
\right\}
\end{eqnarray*}
by the definition of $ \left\{ \mathcal{R}_{n}^{(N)} \right\}_{n=0}^{ \infty } $.
We note $\mathbb{P}_{\mathcal{R}}^{(N)} = \mu_{N} \times P $, the above relation implies
\begin{eqnarray}
&&\mathbb{P}_{\mathcal{R}}^{(N)} \left( \bigcup_{n=1}^{N}\left\{ \mathcal{R}_{n}^{(N)} = ( y_{n} \mod N) \right\} \right)\nonumber\\
&&= \sum_{l=0}^{ N-1} \sum_{n=1}^{N} \frac{1}{N} P \left\{ 
\begin{array}{ll}
l+S_{i} \notin [y_{i}]_{N}, & 1 \leq i \leq n-1, \\
l+S_{n} \in [y_{n}]_{N} & \\
\end{array}
\right\}.\nonumber\\
\label{eq:PNR7}
\end{eqnarray}

For $ l \in \{ 0,1, \ldots , N-1 \}$ and $ n \in \{ 2,3, \ldots ,N \} $, 
we decompose the event $ \{ l+ S_{n} \in [y_{n}]_{N} \}$ according to the 
value of the first hitting time for $[y_{1}]_{N} , [y_{2}]_{N}, 
\ldots , [y_{n}]_{N}$ and the hitting place to obtain 
\begin{eqnarray*}
&&P \{  l+ S_{n} \in [y_{n}]_{N} \}\\
&&\quad = \sum_{j=1}^{n} \sum_{ m \in {\mathbb{Z}} } P \left\{ 
\begin{array}{ll}
l+ S_{i} \notin [y_{i}]_{N} , & \  1 \leq i \leq j-1, \\
\ l+ S_{j}= y_{j}+mN, & \\
 y_{j}+mN+ X_{j+1} + & \cdots + X_{n} \in [y_{n}]_{N} \\
\end{array}
\right\} .
\end{eqnarray*}

The probability in the double summation on the right-hand side above is equal to 
\begin{eqnarray*}
&&P \left\{ 
\begin{array}{ll}
l+ S_{i} \notin [y_{i}]_{N} , & \  1 \leq i \leq j-1, \\
\ l+ S_{j}= y_{j}+mN, & \\
\end{array} 
\right\}\\
&&\quad\times P \left\{ y_{j}+mN+S_{n-j} \in [y_{n}]_{N} \right\}
\end{eqnarray*}
by the Markov property.
It is easy to verify  
that for any $ m \in {\mathbb{Z}}$,
\begin{eqnarray*}
&&P \left\{ y_{j}+mN+S_{n-j} \in [y_{n}]_{N} \right\}\\
&&= P \left\{ S_{n-j} \in [y_{n}-y_{j}]_{N} \right\} \leq p_{n-j}^{(N)}
\end{eqnarray*}
by $ \vert y_{n} - y_{j} \vert \leq n-j $.
Therefore
\begin{eqnarray}
&&P \left\{l+ S_{n} \in [y_{n}]_{N} \right\}\nonumber\\
&&\leq  \sum_{j=1}^{n} P \left\{ 
\begin{array}{ll}
l+ S_{i} \notin [y_{i}]_{N} , & \  1 \leq i \leq j-1, \\
\ l+ S_{j}= [ y_{j}]_{N} & \\
\end{array}
\right\}  p_{n-j}^{(N)} ,\nonumber\\
\label{eq:PNR8}
\end{eqnarray}
for $ l \in \{ 0,1, \ldots , N-1 \} $ and $ n \in \{ 1,2, \ldots ,N \}$. 
By multiplying (\ref{eq:PNR8}) by $1/N$ and summing $(l,n) $ over $ \{ 0,1, \ldots ,N-1 \} \times \{ 1, 2, \ldots ,N \} $, we have
\begin{eqnarray}
&&\sum_{l=0}^{N-1} \sum_{n=1}^{N} \frac{1}{N} P \left\{ l+ S_{n} \in [y_{n}]_{N} \right\}\nonumber\\
&&\leq  \sum_{l=0}^{N-1} \sum_{j=1}^{N} \frac{1}{N} P \left\{ 
 \begin{array}{ll}
 l+ S_{i} \notin [y_{i}]_{N} , & \  1 \leq i \leq j-1, \\
 l+ S_{j}= [ y_{j}]_{N} & \\
\end{array}
\right\}\nonumber\\
&&\quad\times \left( \sum_{i=0}^{N-j} p_{i}^{(N)} \right)\nonumber\\
&&\leq   \mathbb{P}_{\mathcal{R}}^{(N)} \left( \bigcup_{n=1}^{N}\left\{ \mathcal{R}_{n}^{(N)} = ( y_{n} \mod N) \right\} \right)\left( \sum_{i=0}^{N-1} p_{i}^{(N)} \right) . \label{eq:PNR9}
\end{eqnarray}
Here we used (\ref{eq:PNR7}).

By $ \sum_{l=0}^{N-1} P \{ l+S_{n} \in [y]_{N} \}= P \{ S_{n} \in {\mathbb{Z}} \} = 1 \ (n \in {\mathbb{N}}, \ y \in {\mathbb{Z}} )$, 
\begin{equation}
 \sum_{l=0}^{N-1} \sum_{n=1}^{N} \frac{1}{N} P \{ l+ S_{n} \in [y_{n}]_{N} \} 
=1 . \label{eq:PNR10}
\end{equation}

(\ref{eq:PNR9}) and (\ref{eq:PNR10}) imply 
\begin{equation}
1 \leq  \mathbb{P}_{\mathcal{R}}^{(N)} \left( \bigcup_{n=1}^{N}\left\{ \mathcal{R}_{n}^{(N)} = ( y_{n} \mod N) \right\} \right)\left( \sum_{i=0}^{N-1} p_{i}^{(N)} \right) \label{eq:PNR11}
\end{equation}
that is the first inequality in (\ref{eq:PNR}).

For the last inequality in (\ref{eq:PNR}), let $ y_{N+j}=y_{N} \ (j=1,2, \ldots ,N)$. 
The same argument as showing (\ref{eq:PNR11}) (we use $ q_{i}^{(N)} $ instead of $ p_{i}^{(N)} $) gives 
\begin{eqnarray*}
2 &=&  \sum_{l=0}^{N-1} \sum_{n=1}^{2N} \frac{1}{N} P \{ l+ S_{n} \in [y_{n}]_{N} \}\\
&\geq & \mathbb{P}_{\mathcal{R}}^{(N)} \left( \bigcup_{n=1}^{N}\left\{ \mathcal{R}_{n}^{(N)} = ( y_{n} \mod N) \right\} \right)\left( \sum_{i=0}^{N-1} q_{i}^{(N)} \right).
\end{eqnarray*}
\end{Proof}

\begin{corollary}\label{bound;cor2}
For $ N \in {\mathbb{N}} \setminus \{ 1 \} $, 
\begin{eqnarray}
&&\frac{1}{ 1+ \sum_{i=1}^{N-1} P \{ S_{i} \in [0]_{N} \} } \leq \mathbb{P}_{\mathcal{R}}^{(N)} \left( \bigcup_{n=1}^{N}\left\{ \mathcal{R}_{n}^{(N)} = 0 \right\} \right)\nonumber\\
&&\leq \frac{2}{ 1+ \sum_{i=1}^{N-1} P \{ S_{i} \in [0]_{N} \} } .\label{eq:PNR20}
\end{eqnarray}
\end{corollary}

\begin{Proof}
Put $ y_{1}=y_{2}= \cdots =y_{2N}=0 $ in the proof of 
Proposition~\ref{bound;thm2}. Then the same argument as showing 
(\ref{eq:PNR}) gives (\ref{eq:PNR20}). 
\end{Proof}

\begin{corollary}\label{bound;cor3}
For $ N \in {\mathbb{N}} \setminus \{ 1 \} $,
\begin{eqnarray}
&&\frac{1}{ 1+ \sum_{i=1}^{N-1} P \{ S_{i} \in [i]_{N} \} }\nonumber\\
&&\leq \mathbb{P}_{\mathcal{R}}^{(N)} \left( \bigcup_{n=1}^{N}\left\{ \mathcal{R}_{n}^{(N)} = ( n \mod N )  \right\} \right)\nonumber\\
&&\leq \frac{2}{ 1+ \sum_{i=1}^{N-1} P \{ S_{i} \in [i]_{N} \} }.\label{eq:PNR21}
\end{eqnarray}
\end{corollary}

\begin{Proof}
Put $ y_{j}= j \ (j=1,2, \ldots , 2N)$ in the proof of 
Proposition~\ref{bound;thm2}. Then the same argument as showing 
(\ref{eq:PNR}) gives (\ref{eq:PNR21}). 
\end{Proof}

\begin{remark}\label{remark5}
By the same argument as showing (\ref{eq:PNR20}), 
we obtain that for 
$ \tilde{ \epsilon } >0 $ and $ N \geq 1/ \tilde{ \epsilon } $, 
$$ \mathbb{P}_{\mathcal{R}}^{(N)} \left( \bigcup_{n=1}^{N}\left\{ \mathcal{R}_{n}^{(N)} = 0 \right\} \right) 
\leq \frac{ 1+ \tilde{ \epsilon } }{ 1+ \sum_{i=1}^{ \tilde{ \epsilon} N} 
P \{ S_{i} \in [0]_{N} \} } . $$
\end{remark}

\section{Fourier transform}

In this section, we introduce some results concerning one-dimensional random walk. 
\vspace{.1cm}

\begin{proposition}\label{fourier;prop2}
If a one-dimensional random walk satisfies $(A1)$ and $(A3)$, then 
there exist $C_{1}>0$ and $N_{1} \in {\mathbb{N}}$ such that for $n \geq N_{1}$, 
\begin{eqnarray*}
&&\sup_{l \in {\mathbb{Z}}} \left\vert n^{1/ \beta } P \{ S_{n}=l \} - \frac{1}{2 \pi} \int_{ - \infty}^{ + \infty} e^{-c_{*} \vert x  \vert^{ \beta } } \exp \left( -i \frac{  xl }{n^{1/ \beta }} \right) \ dx   \right\vert \\
&&\leq C_{1} n^{ - \delta }, 
\end{eqnarray*}
where $ \delta = \min \{  \varepsilon /(2 \beta ) ,1/2 \}.$
\end{proposition}

\begin{Proof} 
Proposition \ref{fourier;prop2} can be proved by the same procedure as in 
Theorem 1.2.1 of \cite{lawler1991}. 

The Fourier inversion formula for $ \phi^{n}( \theta )$ is
\begin{equation}
n^{1/ \beta} P \{ S_{n}=l \}= \frac{n^{1/ \beta}}{2 \pi} \int_{ - \pi }^{ \pi } 
\phi^{n}( \theta ) e^{-i \theta l}  \ d \theta .
\label{eq:invf}
\end{equation}
By $(A3)$, there exist $C_{*}>0$ and $ r \in (0, \pi )$ such 
that for $ \vert \theta \vert <r$, 
\begin{equation}
 \vert \phi ( \theta ) -(1- c_{*} \vert \theta \vert^{\beta} )
\vert \leq C_{*} \vert \theta \vert^{ \beta + \varepsilon }
\label{eq:a1}
\end{equation} 
and 
\begin{equation}
\vert \phi ( \theta ) \vert \leq 1- \frac{c_{*}}{2} 
\vert \theta \vert^{ \beta} . \label{eq:a2}
\end{equation} 
With $r$, we decompose the right-hand side of 
(\ref{eq:invf}) to obtain 
$$ n^{1/ \beta} P \{ S_{n}=l \} = I(n,l)+ 
J(n,l), $$
where 
\begin{eqnarray*}
&& I(n,l)=  \frac{n^{1/ \beta}}{2 \pi} 
\int_{ \vert \theta \vert <r } \phi^{n}( \theta )
e^{-i \theta l}  \ d \theta, \\
&& J(n,l)= 
 \frac{n^{1/ \beta}}{2 \pi} 
 \int_{ r \leq \vert \theta \vert \leq \pi } \phi^{n}( \theta )
e^{-i \theta l} \ d \theta .
\end{eqnarray*}

A strongly aperiodic random walk $(A1)$ has the property that 
$ \vert \phi ( \theta ) \vert =1$ only when $ \theta $ is a multiple of $2 \pi$ (see \S7 Proposition 8 of \cite{spitzer1976}).
By the definition of $ \phi (\theta )$, $ \vert \phi (\theta ) \vert $ is a continuous function on the bounded closed set $ [-\pi,-r] \cup [r, \pi] $, and $   \vert \phi (\theta ) \vert \leq 1 \ ( \theta \in [ -\pi, \pi] ) $.
Hence, there exists a $ \rho <1$, depending on $r \in (0, \pi] $,  such that 
\begin{equation}
\max_{ r \leq \vert \theta \vert \leq \pi }
\vert \phi (\theta ) \vert \leq \rho. \label{eq:mm}
\end{equation}
By using the above inequality, 
$$ \vert J(n,l) \vert \leq \frac{n^{1/ \beta }}{ 2 \pi} 
\int_{ r \leq \vert \theta \vert \leq \pi }
\vert \phi (\theta )\vert^{n} \ d \theta \leq n^{1/ \beta } 
\rho^{n} . $$

We perform the change of variables 
$ \theta =x / n^{1/ \beta }$, so that 
$$ I(n,l)= \frac{ 1}{2 \pi } \int_{ \vert x \vert <r n^{1/ \beta }}
 \phi^{n} \left( \frac{x}{ n^{1/ \beta }} \right) 
 \exp \left( -i \frac{  xl }{n^{1/ \beta }} \right) \ dx . $$
Put  
$$  \gamma = \min \left\{ \frac{ \varepsilon }
{ 2 \beta ( \beta + \varepsilon +1)} , \ 
\frac{1}{ 2( 2\beta +1 )} \right\}.$$
We decompose $I(n,l)$ as follows: 
\begin{eqnarray*}
I(n,l)&=& \frac{1}{2 \pi} \int_{ - \infty}^{ + \infty} 
e^{-c_{*} \vert x  \vert^{ \beta } }
\exp \left( -i \frac{  xl }{n^{1/ \beta }} \right) \ dx\\
 && + I_{1}(n,l)+I_{2}(n,l)+I_{3}(n,l), 
\end{eqnarray*}

where 
$$
I_{1}(n,l)= \frac{1}{2 \pi} 
\int_{ \vert x \vert \leq n^{\gamma} } 
 \left \{ \phi^{n} \left( \frac{x}{ n^{1/ \beta }} \right) - e^{-c_{*} \vert x \vert^{ \beta } } 
 \right\}  \quad \quad \quad \quad \quad  $$
 $$ \times \exp \left( -i \frac{  xl }{n^{1/ \beta }} \right) \ dx , $$
$$ I_{2}(n,l)= - \frac{1}{2 \pi} \int_{n^{ \gamma } < \vert 
x \vert }  e^{-c_{*} \vert x  \vert^{ \beta } }
\exp \left( -i \frac{  xl }{n^{1/ \beta }} \right) \ dx \quad \quad 
\quad $$
and
$$  I_{3}(n,l)= \frac{ 1}{2 \pi } 
\int_{ n^{ \gamma } <  \vert x \vert <r n^{1/ \beta }}
 \phi^{n} \left( \frac{x}{ n^{1/ \beta }} \right) 
 \exp \left( -i \frac{  xl }{n^{1/ \beta }} \right) \ dx .
$$
Therefore,
$$ \left\vert  n^{1/ \beta } P \{ S_{n}=l \} - 
\frac{1}{2 \pi} \int_{ - \infty}^{ \infty} e^{-c_{*} \vert x  \vert^{ \beta } }
\exp \left( -i \frac{  xl }{n^{1/ \beta }} \right) \ dx   \right\vert  $$
$$ \leq \vert J(n,l) \vert + \sum_{k=1}^{3} \vert I_{k}(n,l) \vert .  \quad 
\quad \quad \quad \quad \quad \quad \quad \quad \quad $$

The proof of Proposition~\ref{fourier;prop2} will be complete if we show that each term in the right-hand side of the above inequality is bounded by a constant (independent of $l$)  multiple 
of $n^{- \delta}$.

If $n$ is large enough, then
the bound $ \vert J(n,l) \vert \leq n^{1/ \beta } \rho^{n} $, which has already been shown above, yields  
$$ \vert J(n,l) \vert \leq n^{- \delta}. $$

With the help of 
\begin{eqnarray}
\vert a^{n}-b^{n} \vert &=& \vert a-b \vert  \left\vert \sum_{j=0}^{n-1} a^{n-1-j}b^{j}\right\vert\nonumber\\
&\leq & n \vert a-b \vert \quad (a,b \in [ -1,1]) \label{eq:anbn}
\end{eqnarray}
and $ \vert \phi( \theta ) \vert \leq 1 \ ( \theta \in [ -\pi ,\pi ])$, 
(\ref{eq:a1}) implies that for $ \vert x \vert < r n^{1/ \beta }$,
\begin{eqnarray*}
&&\left\vert \phi^{n} \left( \frac{x}{ n^{1/ \beta }} \right) - e^{-c_{*} \vert x \vert^{ \beta } } \right\vert \leq n \left\vert  \phi \left( \frac{x}{ n^{1/ \beta }} \right) - e^{-c_{*} \vert x \vert^{ \beta } /n} \right\vert\\
&&\leq n \left\vert  \phi \left( \frac{x}{ n^{1/ \beta }} \right) - \left(1-c_{*} \frac{ \vert x \vert^{ \beta}}{ n} \right) \right\vert\\
&&\hspace{1cm} + n \left\vert \left(1-c_{*} \frac{ \vert x \vert^{ \beta}}{ n} \right) - e^{-c_{*} \vert x \vert^{ \beta } /n} \right\vert\\
&&\leq C_{*} \vert x \vert^{ \beta + \varepsilon}n^{- \varepsilon / \beta } + \frac{c_{*}^{2}}{2} \vert x \vert^{2 \beta} n^{-1}.
\end{eqnarray*}
Thus
\begin{eqnarray*}
\vert I_{1}(n,l) \vert &\leq & \frac{1}{2 \pi } \int_{ \vert x \vert \leq n^{ \gamma } }\left\vert \phi^{n} \left( \frac{x}{ n^{1/ \beta }} \right) - e^{-c_{*} \vert x \vert^{ \beta } } \right\vert \ d \theta\\
&\leq & \frac{1}{ \pi}\left( \frac{C_{*}}{ \beta + \varepsilon +1} + \frac{c_{*}^{2}}{2(2 \beta +1)}\right) n^{- \delta } .
\end{eqnarray*}
It is easy to verify that for $ \vert x \vert < r n^{1/ \beta } $, 
$$ \left\vert \phi^{n} \left( \frac{x}{n^{1/ \beta}} \right) \right\vert \leq 
\left( 1 -\frac{c_{*}}{2} \frac{ \vert x \vert^{ \beta }}{n} \right)^{n} 
\leq e^{ -c_{*} \vert x \vert^{ \beta } /2} $$
by (\ref{eq:a2}), and we obtain that
\begin{eqnarray}
\vert I_{3}(n,l) \vert &\leq & \frac{1}{ 2 \pi } \int_{n^{ \gamma} < \vert x \vert < r n^{1/ \beta }} \left\vert \phi^{n} \left( \frac{x}{n^{1/ \beta}} \right) \right\vert \ dx\nonumber\\
&\leq & \frac{1}{ 2 \pi } \int_{n^{ \gamma} < \vert x \vert } e^{ -c_{*} \vert x \vert^{ \beta } /2} \ dx . \label{eq:i31}
\end{eqnarray}
Moreover, if $n$ is large enough, then 
$$ e^{-c_{*} \vert x \vert^{ \beta }/2} \leq \frac{ 2^{s}}{ c_{*}^{s}}
\vert x \vert^{-s \beta} \quad ( \vert x \vert > n^{ \gamma }), $$
where $ s= (1/ \beta )(1+ 1/(2 \gamma )).$
By replacing the integrand in the right-hand side of the last inequality of (\ref{eq:i31}) with the right-hand side of the above inequality, we obtain
\begin{equation}
 \vert I_{3}(n,l) \vert  \leq \frac{2^{s+1} \gamma}{ \pi c_{*}^{s}}n^{- 1/2 } 
 \leq \frac{2^{s+1} \gamma}{ \pi c_{*}^{s}} n^{- \delta } . 
 \label{eq:i3}
\end{equation}

The same argument as showing (\ref{eq:i3}) gives 
$$ \vert I_{2}(n,l) \vert \leq \frac{1}{ 2 \pi} 
\int_{n^{ \gamma } \leq \vert \theta \vert } e^{-c_{*} \vert x \vert ^{\beta}} \ 
dx \leq \frac{2^{s+1} \gamma}{ \pi c_{*}^{s}}n^{- \delta } .$$
\end{Proof}

Let 
$$ I_{0}( n,l: \beta ,c_{*})= \frac{1}{2 \pi} \int_{ - \infty }^{ + \infty } 
 e^{-c_{*} \vert x  \vert^{ \beta } }
\exp \left( -i \frac{  xl }{n^{1/ \beta }} \right) \ dx $$
appearing in Proposition~\ref{fourier;prop2}. 

\vspace{0.1cm}
\begin{remark}
When a one-dimensional random walk is the strongly aperiodic $(A1)$ 
with $ E[X_{1}]=0$ and 
$E [ \vert X_{1} \vert^{2+ \varepsilon }] < \infty $ for some $ \varepsilon \in 
(0,1) $, it is verified that
$$ \phi ( \theta )= 1 - \frac{E[X_{1}^{2}]}{2} \vert \theta \vert^{2}
+ O( \vert \theta \vert^{2+ \varepsilon }). $$
In this case, $I_{0}(n,l:2,E[X_{1}^{2}]/2)$ can be computed and 
Proposition~\ref{fourier;prop2} gives the following.

{\bf (Local Central Limit Theorem)} 
There exist $\tilde{C}_{1}>0$ and $ \tilde{N}_{1} \in {\mathbb{N}}$ such that for $n \geq \tilde{N}_{1}$,
$$ 
\sup_{l \in {\mathbb{Z}}} \left\vert n^{1/2 } P \{ S_{n}=l \} - \frac{1}{  \sqrt{ 2 E[X_{1}^{2}] \pi}} \exp \left( - \frac{l^{2}}{2 E[X_{1}^{2}]n} \right) \right\vert 
$$
\begin{equation}
\leq \tilde{C}_{1} n^{ - \delta },  \quad \quad \quad \quad \quad \quad 
 \quad \quad \quad \quad \quad \quad \quad \quad \quad \quad \quad \quad 
 \label{eq:LCLT}
\end{equation}
where $ \delta = \min \{  \varepsilon /4 ,1/2 \}.$
(See Remark after Proposition 7.9 in \cite{spitzer1976}.)
\vspace{0.1cm}

It is easy to see
$$ I_{0}( n,l: 1,c_{*} )= \frac{1}{ \pi } \frac{ c_{*}}{ c_{*}^{2}+ (l/n)^{2}} 
\quad ( n \in {\mathbb{N}}, l \in {\mathbb{Z}}, c_{*} >0) $$
and we have the following corollary of Proposition~\ref{fourier;prop2}. 
\end{remark}

\begin{corollary}\label{fourier;cor4}
If a one-dimensional random walk satisfies $(A1)$ and $(A3)$ with $ \beta =1$, then 
there exist $C_{2}>0$ and $N_{2} \in {\mathbb{N}}$ such that for $n \geq N_{2}$, 
$$ \sup_{l \in {\mathbb{Z}}} \left\vert n P \{ S_{n}=l \} - 
\frac{1}{ \pi } \frac{ c_{*}}{ c_{*}^{2}+ (l/n)^{2}}  \right\vert 
\leq C_{2} n^{ - \delta }, $$
where $ \delta = \min \{  \varepsilon /2 ,1/2 \}.$
\end{corollary}

We perform the change of variables $ t= c_{*} x^{ \beta } $  , so that 
$$ I_{0}(n,0: \beta ,c_{*} ) = \frac{1}{\pi} \int_{0}^{+ \infty } 
e^{- c_{*} x^{ \beta }} \ dx = \frac{1}{ \beta c_{*}^{1 / \beta } \pi } \Gamma \left( \frac{1}{ \beta } \right) .
$$
With the help of the above calculation, Proposition~\ref{fourier;prop2} gives the following corollary.

\begin{corollary}\label{fourier;cor5}
If a one-dimensional random walk satisfies $(A1)$ and $(A3)$, then 
there exist $C_{3}>0$ and $N_{3} \in {\mathbb{N}}$ such that for $n \geq N_{3}$, 
$$ \left\vert n^{1 / \beta} P \{ S_{n}=0 \} - 
\frac{1}{ \beta c_{*}^{1 / \beta } \pi } \Gamma \left( \frac{1}{ \beta } \right)   \right\vert 
\leq C_{3} n^{ - \delta }, $$
where $ \delta = \min \{  \varepsilon /2 \beta ,1/2 \}.$
\end{corollary}

\begin{proposition}\label{prop:one-dimensional}
If a one-dimensional random walk satisfies $(A2)$, then for 
$ l \in {\mathbb{Z}} $ and $ n \in  \{ 0 \}  \cup {\mathbb{N}} $,
\begin{eqnarray}
&&P \left\{ S_{n} \in [l]_{N} \right\}\nonumber\\
&&= \frac{1}{N}+ \frac{2}{N} \sum_{1 \leq j \leq( N-1)/2 }
\phi^{n}\left( \frac{2j \pi}{N} \right) \cos \left( \frac{2j \pi}{N} l \right) + J_{N}(n,l),\nonumber\\
\label{eq:876}
\end{eqnarray}
where 

\begin{eqnarray*}
J_{N}(n,l)= \left\{
\begin{array}{ll}
(1/N)  \phi^{n}(\pi)\cos( \pi l ) & \quad (\mbox{ if }  N \mbox{ is even }) \\
0                                 & \quad ( \mbox{ if } N \mbox{ is odd }).\\
\end{array}
\right.
\end{eqnarray*}


\end{proposition}

\begin{Proof}
By the definition of $\phi (\theta)$,
\begin{eqnarray*}
\phi ^{n}(\theta) = \sum_{k \in \mathbb{Z}}e^{i\theta k}P\left\{S_{n} = k\right\}.
\end{eqnarray*}
Thus
\begin{eqnarray*}
&&\phi ^{n}\left(\frac{2j\pi}{N} \right) = \sum_{k \in \mathbb{Z}}e^{2ij\pi k/N}P\left\{S_{n} = k\right\}\\
&&= \sum_{\tilde{l} = 0}^{N-1}\sum_{m \in \mathbb{Z}}e^{2ij\pi (\tilde{l}+mN)/N}P\left\{S_{n} = \tilde{l}+mN\right\}\\
&&= \sum_{\tilde{l} = 0}^{N-1}e^{2ij\pi \tilde{l}/N}P\left\{S_{n} \in [\tilde{l}]_{N}\right\}.
\end{eqnarray*}
Then,
\begin{eqnarray*}
&&\sum_{j=0}^{N-1}e^{-2ij\pi l/N}\phi^{n}\left(\frac{2j\pi}{N}\right) = \sum_{\tilde{l} = 0}^{N-1}\sum_{j=0}^{N-1}e^{2ij\pi (\tilde{l}-l)/N}P\left\{S_{n} \in [\tilde{l}]_{N}\right\}\\
&&= NP\left\{S_n \in [l]_{N} \right\}
\end{eqnarray*}
since
\begin{eqnarray*}
\sum_{j=0}^{N-1} e^{2ij\pi (\tilde{l}-l)/N} = \left\{
\begin{array}{ll}
 N &\quad \tilde{l} = l\\
 0 &\quad \tilde{l} \neq l 
\end{array}
\right. .
\end{eqnarray*}
Therefore,
\begin{eqnarray*}
P\left\{S_{n} \in [l]_{N}\right\} &=& \frac{1}{N}\sum_{j = 0}^{N-1}\phi ^{n}\left(\frac{2j\pi}{N}\right)e^{-2j\pi il/N}\\
&=& \frac{1}{N}\sum_{j = 0}^{N-1}\phi ^{n}\left(\frac{2j\pi}{N}\right)\cos \left(\frac{2j\pi l}{N}\right).
\end{eqnarray*}
We note that $\phi ^{n}(\theta) \in \mathbb{R}$ and 
\begin{eqnarray*}
\frac{1}{N}\sum_{j = 0}^{N-1}\phi ^{n}\left(\frac{2j\pi}{N}\right)\cos \left(\frac{2j\pi l}{N}\right) \in \mathbb{R}
\end{eqnarray*}
by $(A2)$.
So we have
\begin{eqnarray}
&&\phi ^{n}\left(\frac{2m\pi}{N}\right)\cos \left(\frac{2m\pi l}{N}\right)\nonumber\\
&&= \phi ^{n}\left(\frac{2(N-m)\pi}{N}\right)\cos \left(\frac{2(N-m)\pi l}{N}\right).\label{eq:yi01}
\end{eqnarray}
Let $N$ be an even number.
Then, by (\ref{eq:yi01}),
\begin{eqnarray*}
&&P\left\{S_{n} \in [l]_{N}\right\}\\
&&= \frac{1}{N}\phi ^{n}\left(0\right)\cos \left(0\right)\\
&&\quad + \frac{2}{N}\sum_{1 \leq j \leq (N-1)/2}\phi ^{n}\left(\frac{2j\pi}{N}\right)\cos \left(\frac{2j\pi l}{N}\right)\\
&&\quad+ \frac{1}{N}\phi ^{n}\left(\pi\right)\cos \left(\pi l\right)\\
&&= \frac{1}{N} + \frac{2}{N}\sum_{1 \leq j \leq (N-1)/2}\phi ^{n}\left(\frac{2j\pi}{N}\right)\cos \left(\frac{2j\pi l}{N}\right)\\
&&\quad+ \frac{1}{N}\phi ^{n}\left(\pi\right)\cos \left(\pi l\right).
\end{eqnarray*}
Therefore, we have (\ref{eq:876}) for every even number $N$. 
The proof of (\ref{eq:876}) for odd number is similar and is omitted. 
\end{Proof}

\section{Proof of Theorem~\ref{main_theorem}}

In this section we prove Theorem~\ref{main_theorem}. To prove it, we introduce the 
following Proposition. 

\begin{proposition}\label{proof;prop4}
Assume $(A1)$, $(A2)$ and $(A3)$.

If $ \beta \in (0,1)$, then there exists a constant 
$c_{8}>0$ such that
\begin{equation}
 \sum_{i=0}^{ N-1} p_{i}^{(N)} \leq c_{8}. \label{eq:pop41}
\end{equation}

If $ \beta =1 $, then there exists a constant 
$c_{9} >0$ such that
\begin{equation}
\sum_{i=0}^{ N-1} p_{i}^{(N)} \leq \frac{1}{c_{*} \pi } 
\log N + c_{9}. \label{eq:pop42}
\end{equation}

If $ \beta \in (1,2]$, then there exists a constant 
$c_{10}>0$ such that
\begin{equation}
 \sum_{i=0}^{ N-1} p_{i}^{(N)} \leq c_{10} 
 N^{( \beta -1) / \beta}. \label{eq:pop43}
\end{equation}
\end{proposition}

\begin{Proof} 

There exist $C_{*}$ and $ r \in (0 ,\pi )$ such that 
for $ \vert \theta \vert < r$, 
\begin{equation}
 \vert \phi ( \theta ) -(1- c_{*} \vert \theta \vert^{ \beta }) \vert 
\leq C_{*} \vert \theta \vert^{ \beta + \varepsilon } \label{eq:VPT}
\end{equation} 
by $(A3)$.
We can choose $ r_{*} \in (0,r] $ small enough so that 
\begin{eqnarray}
C_{*} r_{*} ^{ \varepsilon } \leq  \frac{1}{2} c_{*} \ \mbox{ and } 
\  c_{*} r_{*}^{ \beta } \leq \frac{1}{3}. \label{eq:C*r*}
\end{eqnarray}
Then for $ \vert \theta \vert <r_{*} $, 
\begin{equation}
\frac{1}{2} c_{*} \vert \theta \vert^{ \beta }
 \leq \vert 1 -\phi (\theta ) \vert  
 \label{eq:V1}
\end{equation}
and 
\begin{equation}
 \vert 1 -\phi (\theta ) \vert   \leq  
\frac{3}{2}c_{*} \vert \theta \vert^{ \beta } 
 \leq \frac{1}{2}. \label{eq:V2}
\end{equation} 
There exists a $ \rho_{*} \in [0,1) $, depending on $r_{*}$, such that 
\begin{equation}
\max_{ r_{*} \leq \vert \theta \vert \leq \pi }
\vert \phi (\theta ) \vert \leq \rho_{*} \label{eq:mmm}
\end{equation}
by the same reason as (\ref{eq:mm}).
(Here we used the condition $(A1)$.)

Using Proposition~\ref{prop:one-dimensional} and (\ref{eq:mmm}), we obtain that for 
$ i \in \{ 1, 2, \ldots , N-1 \}$,
\begin{eqnarray*}
p_{i}^{(N)} &=& \max_{\vert l \vert \leq i }P\left\{S_{i} \in [l]_{N} \right\}\\
&\leq & \frac{1}{N} +   \sum_{1 \leq j \leq (N-1)/2}  \frac{2}{N} \left\vert\phi\left(\frac{2j\pi}{N}\right)\right\vert^{i} + \vert J_{N}(i,0)\vert\\
&\leq & \frac{1}{N} + \sum_{1 \leq j < (r_{*}/(2 \pi))N} \frac{2}{N} \left\vert\phi\left(\frac{2j\pi}{N}\right)\right\vert^{i} + \rho_{*}^{i}.
\end{eqnarray*}
Therefore
\begin{equation}
 \sum_{i=0}^{ N-1} p_{i}^{(N)} \leq 1+ \Phi_{N}
+ \frac{1}{ 1- \rho_{*}} , \label{eq:spn1}
\end{equation}
where 
$$ \Phi_{N}=  \sum_{1 \leq j <  (r_{*}/(2 \pi))N } 
\frac{2}{N} 
\frac{ 1-\left\vert\phi\left(\frac{2j\pi}{N}\right)\right\vert^{N}}{1-
\left\vert \phi\left(\frac{2j\pi}{N}\right) \right\vert }. $$

Because of $(A2)$, $ \phi (\theta )$ takes a real number. Then  (\ref{eq:V1}), 
(\ref{eq:V2}) and $(A1)$ mean that 
\begin{equation}
\frac{1}{2} < \phi ( \theta ) = \vert 
 \phi ( \theta )  \vert <1 \quad (\theta \in (-r_{*} ,0) \cup 
(0, r_{*}) ) \label{eq:V3}
\end{equation}
and 
\begin{equation}
\Phi_{N}  \leq \sum_{
1 \leq j <  (r_{*}/(2 \pi))N }  \frac{2}{N} \frac{ 1 }{1-
\phi\left(\frac{2j\pi}{N}\right) }. \label{eq:spn2}
\end{equation}

We will calculate $ \Phi_{N}$ in the case $ \beta \in (0,1]$.  
By (\ref{eq:spn2}), we decompose the right-hand side of the above to obtain 
\begin{equation}
 \sum_{
1 \leq j <  (r_{*}/(2 \pi))N } \frac{2}{N} \frac{ 1 }{1-
\phi\left(\frac{2j\pi}{N}\right) } = \tilde{\Phi}_{N}+ E_{N} , 
\label{eq:spn3}
\end{equation}
where 
$$  \tilde{\Phi}_{N}= 
\frac{2^{1- \beta }}{ \pi^{ \beta }c_{*}} N^{ \beta -1} \sum_{
1 \leq j <  (r_{*}/(2 \pi))N } j^{ - \beta } , \quad \quad \quad \quad 
\quad \quad \quad $$
$$ E_{N}=  \sum_{
1 \leq j <  (r_{*}/(2 \pi))N }  \frac{2}{N} \left(
 \frac{ 1 }{1-\phi\left(\frac{2j\pi}{N}\right) }
- \frac{ 1 }{c_{*}
\left(\frac{2j\pi}{N}\right)^{ \beta } } \right) .$$

To estimate $E_{N}$, we use (\ref{eq:VPT}) and (\ref{eq:V1}) which imply that for $ j \in [1, (r_{*}/(2 \pi))N) \cap \mathbb{Z}$,
\begin{eqnarray*}
&&\frac{2}{N} \left\vert  \frac{ 1 }{1-\phi\left(\frac{2j\pi}{N}\right) }- \frac{ 1 }{c_{*}\left(\frac{2j\pi}{N}\right)^{ \beta } } \right\vert\\
&&=  \frac{2}{N}\frac{ \left\vert 1-\phi\left(\frac{2j\pi}{N}\right) - c_{*} \left(\frac{2j\pi}{N}\right)^{ \beta } \right\vert }{\left\vert 1-\phi\left(\frac{2j\pi}{N}\right) \right\vert \cdot\left\vert c_{*} \left(\frac{2j\pi}{N}\right)^{ \beta } \right\vert }  \leq c_{11} N^{ \beta - \varepsilon -1} j^{ \varepsilon - \beta } ,
\end{eqnarray*}
 where $c_{11}= 2^{2+ \varepsilon - \beta } \pi^{ \varepsilon - \beta } 
 C_{*}/c_{*}^{2} .$
By noticing that $ 1+ \varepsilon - \beta >0$, 
$$ \sum_{ 1 \leq j <  (r_{*}/(2 \pi))N  } j^{ \varepsilon - \beta } 
\leq \int_{0}^{N} x^{ \varepsilon - \beta } \ dx = 
\frac{ N^{ 1+ \varepsilon - \beta }}{  1+ \varepsilon - \beta } .$$
Thus 
\begin{equation}
\vert E_{N} \vert \leq c_{11}/( 1+ \varepsilon - \beta ). \label{eq:spn4}
\end{equation}

It is easy to see that
\begin{eqnarray}
\tilde{\Phi}_{N} &\leq & \frac{2^{1- \beta }}{ \pi^{ \beta }c_{*}} N^{ \beta -1}\left( 1+ \int_{1}^{N} x^{ - \beta } \ dx \right)\nonumber\\
&\leq & \left\{ 
 \begin{array}{ll}
 \displaystyle{ \frac{2^{1- \beta }}{\pi^{ \beta }c_{*}(1- \beta)} }& ( \beta \in (0,1)) \\
 \displaystyle{\frac{1}{ \pi c_{*}} \log N +\frac{1}{ \pi c_{*}} } & ( \beta =1) .
 \end{array} 
\right.  \label{eq:spn5}
\end{eqnarray}

Put the pieces ((\ref{eq:spn1}), (\ref{eq:spn2})-(\ref{eq:spn5})) together, we have (\ref{eq:pop41}) and (\ref{eq:pop42}).

In the case $ \beta \in (1,2]$, we use (\ref{eq:V3}) to obtain 
\begin{equation}
 \Phi_{N} \leq \Phi_{N}^{(1)} +\Phi_{N}^{(2)}, \label{eq:spn12}
\end{equation}
where $ N( \beta )= \min \{ N^{( \beta -1) / \beta } ,(r_{*}/( 2 \pi))N \} $ 
and 
$$ \Phi_{N}^{(1)}= 
\sum_{ 1 \leq j  < N( \beta)}  \frac{2}{N} \frac{ \left\vert 1- \phi\left(\frac{2j\pi}{N}\right)^{N} \right\vert } { \left\vert 1 - \phi\left(\frac{2j\pi}{N}\right) \right\vert } , \quad \quad \quad $$
$$ \Phi_{N}^{(2)}=  
\sum_{ N( \beta ) \leq j < (r_{*}/( 2 \pi))N } 
\frac{2}{N} \frac{ 1 } { \left\vert 1 - \phi\left(\frac{2j\pi}{N}\right) \right\vert } .  $$
We use (\ref{eq:anbn})(set $n=N$ and $a=1,b= \phi \left(\frac{2j\pi}{N}\right)$), then
\begin{equation}
\Phi_{N}^{(1)} \leq 2N( \beta ) \leq 2 N^{( \beta -1)/ \beta } . 
\label{eq:spn13}
\end{equation}
We notice that $  \beta -1>0$, (\ref{eq:V1}) gives
\begin{eqnarray}
\Phi_{N}^{(2)} &\leq & \frac{2^{2- \beta }}{c_{*} \pi^{ \beta }} N^{ \beta -1} \left( \sum_ {N( \beta ) \leq j < (r_{*}/( 2 \pi))N }  j^{- \beta }\right)\nonumber\\
&\leq & \frac{2^{2- \beta }}{c_{*} \pi^{ \beta }}N^{ \beta -1} \left(N^{- \beta +1} + \int_{ N^{( \beta -1)/ \beta }}^{ + \infty } x^{ - \beta } \ dx \right)\nonumber\\
&\leq & \frac{2^{2- \beta }}{c_{*} \pi^{ \beta }} \left( 1+\frac{1}{ \beta -1} \right) N^{( \beta -1 )/ \beta }.\label{eq:spn14}
\end{eqnarray}
Put the pieces ((\ref{eq:spn1}), (\ref{eq:spn12})-(\ref{eq:spn14})) together, we have (\ref{eq:pop43}).
\end{Proof}

It remains to show the last inequality in (\ref{eq:r2}). To achieve this,  
we will use Proposition~\ref{bound;thm2} and Corollary~\ref{fourier;cor4}.

There exist $C_{2}>0$ and $N_{2} \in \mathbb{N}$ 
such that for $ i \geq N_{2}$ and $ l \in \mathbb{Z}$, 
$$ P \{S_{i} =l \} \geq \frac{1}{ \pi }\frac{c_{*}}{c_{*}^{2}+(l/i)^{2}} \frac{1}{i} 
- C_{2} i^{-1- \delta } $$
by Corollary~\ref{fourier;cor4}.
Let
\begin{eqnarray*}
c_{12}:= \frac{1}{ \pi} \frac{ c_{*}}{c_{*}^{2}+1} \log N_{2} +C_{2} 
\sum_{i=N_{2}}^{\infty } i^{-1- \delta } .
\end{eqnarray*}
We can choose $ N_{*} \in \mathbb{N} $ large enough so that 
\begin{eqnarray*}
\frac{1}{2} \frac{1}{ \pi} \frac{ c_{*}}{c_{*}^{2}+1} \log N_{*} \geq 
c_{12}.
\end{eqnarray*}
Then for $N \geq N_{*}+1$,
\begin{eqnarray}
\sum_{i=0}^{N-1} q_{i}^{(N)} &\geq &  \sum_{i=N_{2}}^{N-1}\min_{ \vert l \vert \leq i} P \{ S_{i} =l \}\nonumber\\
&\geq &  \frac{1}{ \pi }\frac{c_{*}}{c_{*}^{2}+1} \sum_{i=N_{2}}^{N-1}\frac{1}{i} - C_{2} \sum_{i=N_{2}}^{ \infty} i^{ -1- \delta }\nonumber\\
&\geq &  \frac{1}{ \pi }\frac{c_{*}}{c_{*}^{2}+1} \log N -c_{12}\nonumber\\
&\geq & \frac{1}{2} \frac{1}{ \pi }\frac{c_{*}}{c_{*}^{2}+1} \log N. \label{eq:si0N}
\end{eqnarray}
It follows from Proposition~\ref{bound;thm2} and (\ref{eq:si0N}) that 
for $ N \in [N_{*}+1, + \infty ) \cap \mathbb{N} $ and $y_{1},y_{2}, \ldots , y_{N} \in {\mathbb{Z}} $ 
with $ \vert y_{n} -y_{n+1} \vert \leq 1 \ (n=1,2, \ldots , N-1)$, 
$$  \mathbb{P}_{\mathcal{R}}^{(N)} \left( \bigcup_{n=1}^{N}
\{ \mathcal{R}_{n}^{(N)} = ( y_{n} \mod N) \} \right)  \leq 
\frac{ \frac{4 \pi (c_{*}^{2}+1)}{c_{*}}}{ \log N} . $$
It is clear that $ \mathbb{P}_{\mathcal{R}}^{(N)} \left( \bigcup_{n=1}^{N}
\{ \mathcal{R}_{n}^{(N)} = ( y_{n} \mod N) \} \right)  $ is bounded by 1.
Put $ c_{3}= \max \{ 4 \pi (c_{*}^{2}+1)/c_{*},  \log N_{*} \} $. 
The last inequality in (\ref{eq:r2}) holds.

The proof of Theorem~\ref{main_theorem} is complete.

\section{Conclusion and Future works}
We formalized the Hunter vs Rabbit game using the random walk framework.
We generalize a probability distribution of the rabbit's strategy using four assumptions.
We have the general lower bound formula of a probability that the rabbit is caught.
Let $P\left\{ X_{1} = k\right\} = O(k^{-\beta -1})$.
If $\beta \in (0,1)$, the lower bound of a probability that the hunter catches the rabbit is $c_1$ where $c_1 > 0$ is a constant.
If $\beta = 1$, the lower bound of a probability that the rabbit is caught is $\frac{1}{\frac{1}{c_{*}\pi}\log N + c_2}$ where $c_{2}$ and $c_{*}$ are constants defined by the given strategy.
If $\beta \in (1,2]$, the lower bound of a probability that the rabbit is caught is $\frac{c_4}{N^{(\beta-1)/\beta}}$ where $c_4 > 0$ is a constant defined by the given strategy.\par
We show experimental results for three examples of the rabbit strategies.
We can confirm our bounds formula, and asymptotic behavior of those bounds
\begin{eqnarray*}
\lim_{N \rightarrow \infty}\left( \frac{1}{ c_{*} \pi } \log N \right) \mathbb{P}_{\mathcal{R}}^{(N)} \left( \bigcup_{n=1}^{N}\{ \mathcal{R}_{n}^{(N)} = 0  \} \right) = 1.
\end{eqnarray*}
\par
In this paper, we consider the lower bound of a probability that the rabbit is caught to show the worst expected value of time until the rabbit caught.
Our motivation is to find the best strategy of the rabbit.
Our results help to find the best strategy of the rabbit.
On the other hands, what is the best strategy of the hunter?
And what is the worst strategy of the hunter?
Future works include to show the best strategy of the hunter is $Y_{j+1} = Y_{j} + 1$, and the worst strategy of the hunter is $Y_{j} = \mathcal{H}_{0}^{(N)}$ for any $j$.

\section{Acknowledgment}
I would like to express my deepest gratitude to Professor Hiroyuki Ochiai for his valuable advice and guidance.
I would like to thank Mr. Norikazu Ishii for his help.


\section*{Appendix} 

(A) 
{\bf Proof of Proposition \ref{prop;add1}.}
The first inequality in (\ref{eq:00}) comes from 
(\ref{eq:r3}) in Theorem~\ref{main_theorem}. To prove the last inequality in (\ref{eq:00}), 
we will use Corollary~\ref{bound;cor2} and \ref{fourier;cor5} instead of Proposition~\ref{bound;thm2} and Corollary~\ref{fourier;cor4}. The same 
argument as showing the  last inequality in (\ref{eq:r3}) 
gives the last inequality in (\ref{eq:00}). 
\qed

{\bf Proof of Proposition \ref{prop;add2}.}
We consider the case when $X_{1}$ takes three values $ -1,0,1$ with equal probability.
In this case, $X_{1}$ satisfies $(A1)$, $(A2)$ and 
$$
\phi ( \theta )= 1 -\frac{1}{3} \vert \theta \vert^{2}+ O( \vert \theta \vert^{4}).
$$
We can show that there exist $\tilde{C}_{1}>0$ and $ \tilde{N}_{1} \in {\mathbb{N}}$ such that for $i \geq \tilde{N}_{1}$ and $l \in \mathbb{Z} $,
\begin{equation}
 P \{ S_{i} = l \} \leq \frac{ \sqrt{3} }{2 \sqrt{ \pi }} \frac{1}{ i^{1/2}} 
 \exp \left( - \frac{3l^{2}}{4i} \right) + \tilde{C}_{1} i^{-1}
 \label{eq:LCLT1}
\end{equation}
by (\ref{eq:LCLT}).
We notice that $ P \{ \vert X_{1} \vert \leq 1 \}=1$,
then we obtain that for $ N \in \mathbb{N} \setminus \{ 1 \} $, 
$$ 1+ \sum_{i=1}^{N-1} P \{ S_{i} \in [i]_{N} \} \quad 
\quad \quad \quad \quad \quad \quad \quad \quad \quad
 \quad \quad \quad $$
$$ = 
1+ \sum_{i=1}^{ N-1} P \{ S_{i}=i \}  + \sum_{ N/2 \leq i \leq N-1} 
P \{ S_{i}= i-N \} $$
and 
$$ \sum_{i=1}^{ N-1} P \{ S_{i}=i \} = \sum_{i=1}^{ N-1}
\left( \frac{1}{3} \right)^{i} \leq \frac{1}{2} .$$
With the help of $ e^{-x} \leq 1/x \ (x>0)$, 
(\ref{eq:LCLT1}) implies that for $ N \geq 2 \tilde{N}_{1} $,
\begin{eqnarray*}
&&\sum_{N/2 \leq k \leq N-1 } P \{ S_{k}= k-N \}\\
&&\leq \sum_{N/2 \leq k \leq N-1 } \left\{\frac{ \sqrt{3}}{ 2 \sqrt{ \pi} } \frac{1}{k^{1/2}}\exp \left( - \frac{3(k-N)^{2}}{4k} \right) + \tilde{C}_{1} k^{-1}  \right\}\\
&&\leq  \sqrt{ \frac{3}{2 \pi} } \frac{1}{ N^{1/2}}  \sum_{ 1 \leq k \leq N/2}\exp \left( - \frac{3k^{2}}{ 4N} \right) + \tilde{C}_{1}\sum_{ 1 \leq k \leq N/2} \frac{2}{N}\\
&&\leq \sqrt{ \frac{3}{2 \pi} } \frac{1}{ N^{1/2}}\left( \sum_{1 \leq k \leq N^{1/2}} 1 + \sum_{ N^{1/2} < k } \frac{4N}{3k^{2}} \right) + 2 \tilde{C}_{1}\\
&&\leq  \sqrt{ \frac{3}{2 \pi} } + \frac{2 \sqrt{2}}{ \sqrt{3 \pi}} N^{1/2} \left( \frac{1}{N} + \int_{N^{1/2}}^{+ \infty} \frac{1}{x^{2}} \ dx \right)+ 2 \tilde{C}_{1}\\
&&\leq c_{13},
\end{eqnarray*}
where $c_{13}= \sqrt{ 3/ (2 \pi)}+ 4 \sqrt{2}/ \sqrt{3 \pi }+2 \tilde{C}_{1}. $ 
Thus for $ N \in  \mathbb{N} \setminus \{ 1 \}$,  
$$
1+ \sum_{i=1}^{N-1} P \{ S_{i} \in [i]_{N} \} \leq 
\max \{ 2\tilde{N}_{1} , (3/2)+c_{13} \}. $$
Combining the above inequality with Corollary~\ref{bound;cor3}, we have (\ref{eq:000}).
\qed
\vspace{.2cm}

(B) To obtain (\ref{eq:AAAA}), we use the formula 
\begin{equation}
\int_{0}^{ + \infty } \frac{ \sin bx }{ x^{ \alpha }} \ dx =
\frac{ \pi b^{ \alpha -1}}{ 2 \Gamma ( \alpha ) \sin ( \alpha \pi /2)}
\label{eq:BB}
\end{equation}
for $ \alpha \in (0,2) $ and $ b>0$. 
By the definition of $X_{1}$, 
$$ 1 - \phi ( \theta ) = \frac{1}{a} \sum_{k=1}^{ \infty } (1- \cos  \vert 
\theta \vert k ) 
\frac{1}{ k^{ \beta +1 } }. $$
A simple calculation shows that the absolute value of the difference between the 
right-hand side of the above and 
$$ \frac{1}{a} \int_{0}^{ + \infty } \frac{ 1- \cos \vert \theta 
\vert x}{ x^{ \beta +1} } \ dx $$
is bounded by a constant multiple of $ \vert \theta \vert^{ \beta +(2 - \beta )/2}.$ 
It remains to show that 
\begin{equation}
 \frac{1}{a} \int_{0}^{ + \infty } \frac{ 1- \cos \vert \theta 
 \vert x}{ x^{ \beta +1} } \ dx = 
 \frac{ \pi }{ 2a} \frac{ \vert \theta \vert^{ \beta }}{ \Gamma ( \beta +1) 
 \sin ( \beta \pi /2) } . \label{eq:BBB}
\end{equation}
We perform integration by part for the left-hand side of (\ref{eq:BBB}) and 
use (\ref{eq:BB}). Then we have (\ref{eq:BBB}) and (\ref{eq:AAAA}). 

\vspace{.2cm}

(C) {\bf Proof of (\ref{eq:LIM}).} Let $ \epsilon >0 $ be fixed.
By Corollary 4, there exist $C_{2} >0$ and $N_{2} \in \mathbb{N}$ such that 
for $ i \geq N_{2}$, 
\begin{equation}
 P \{ S_{i} =0 \} \geq 
\frac{1}{ c_{*} \pi } \frac{1}{i} - C_{2}i^{ -1- \delta }. \label{eq:ap}
\end{equation}
(\ref{eq:ap}) implies that for $  N \geq  (4/ \epsilon)(N_{2}+1) $,
\begin{eqnarray} 
&&1+ \sum_{1 \leq i \leq ( \epsilon /4)N} P \{ S_{i} \in [0]_{N} \} \geq \sum_{N_{2} \leq i \leq  ( \epsilon /4)N} P \{ S_{i} =0 \}\nonumber\\
&&\geq \sum_{N_{2} \leq i \leq  ( \epsilon /4)N} \left( \frac{1}{ c_{*} \pi } \frac{1}{i} - C_{2}i^{ -1- \delta } \right)\nonumber\\
&&\geq \frac{1}{ c_{*} \pi } \int_{N_{2}}^{ ( \epsilon /4) N} \frac{1}{x} \ dx - C_{2} \left( \frac{1}{N_{2}^{1+ \delta}} + \int_{N_{2}}^{ + \infty } x^{ -1- \delta } \ dx \right)\nonumber\\
&&= \frac{1}{ c_{*} \pi }  \log N + \frac{1}{ c_{*} \pi } \log \epsilon -c_{14},\label{eq:app}
\end{eqnarray}
where $c_{14}= (1/ ( c_{*} \pi )) \log 4 + (1/ ( c_{*} \pi )) \log N_{2} +
C_{2} \{ 1/ N_{2}^{1+ \delta }+  1/ ( \delta N_{2}^{ \delta }) \}.$

We can choose $N_{4} \in \mathbb{N}$ which satisfies  
\begin{equation}
  \min \left\{\frac{1}{2}, \frac{ \epsilon }{ 8} \right\} 
 \frac{1}{ c_{*} \pi } \log N_{4} \geq \left\vert - \frac{1}{ c_{*} \pi}
\log \epsilon + c_{14} \right\vert \label{eq:CON1}
\end{equation}
and 
\begin{equation}
 \frac{ \epsilon }{4} \frac{1}{ c_{*} \pi } \log N_{4} \geq c_{2} , 
 \label{eq:CON2}
\end{equation}
where $c_{2}$ is the same constant in (\ref{eq:r2}).

Combining Remark 5 with (\ref{eq:app}) and using the left-hand side of (\ref{eq:r2}), we obtain that for $N \geq \max \{ N_{4},
 (4/ \epsilon)(N_{2}+1) \}$, 
$$ \frac{1}{ \frac{1}{ c_{*} \pi} \log N +c_{2} } \leq 
 \mathbb{P}_{\mathcal{R}}^{(N)} \left( \bigcup_{n=1}^{N}
\{ \mathcal{R}_{n}^{(N)} =0 \} \right) \quad \quad $$
$$ \quad \quad \quad \quad \quad \quad \quad 
\leq \frac{1+ ( \epsilon /4) }
{\frac{1}{ c_{*} \pi }  \log N + \frac{1}{ c_{*} \pi } \log \epsilon -
c_{14}}. $$
Hence for $N \geq \max \{ N_{4}, (4/ \epsilon)(N_{2}+1) \}$, 
$$ \left\vert \left( \frac{1}{c_{*} \pi } \log N \right) 
 \mathbb{P}_{\mathcal{R}}^{(N)} \left( \bigcup_{n=1}^{N}
\{ \mathcal{R}_{n}^{(N)} =0 \} \right) -1 \right\vert $$
$$ \leq E_{N}^{(1)} + E_{N}^{(2)}, \quad \quad \quad \quad 
\quad \quad \quad \quad \quad \quad \quad $$
where 
$$ E_{N}^{(1)}= \left\vert  
\frac{ \frac{1}{ c_{*} \pi} \log N}{ \frac{1}{ c_{*} \pi} \log N +c_{2} } 
-1 \right\vert \quad \quad \quad \quad \quad $$
and 
$$ E_{N}^{(2)}= \left\vert \frac{(1+ ( \epsilon /4)) \frac{1}{ c_{*} \pi }  \log N}
{ \frac{1}{ c_{*} \pi }  \log N + \frac{1}{ c_{*} \pi } \log \epsilon -
c_{14}} -1 \right\vert . $$

The proof is complete if we show that for \\
$N \geq \max \{ N_{4}, (4/ \epsilon)(N_{2}+1) \}$, 
\begin{equation}
 E_{N}^{(1)} +  E_{N}^{(2)} \leq \epsilon . \label{eq:CON3}
\end{equation}

We use (\ref{eq:CON2}), then 
$$ E_{N}^{(1)} \leq \frac{c_{2}}{ \frac{1}{ c_{*} \pi }  \log N } 
\leq \frac{ \epsilon }{4} 
$$
for $ N \geq \max \{ N_{4}, (4/ \epsilon)(N_{2}+1) \} .$ 
We can show that
\begin{eqnarray*}
E_{N}^{(2)} &\leq & \frac{ ( \epsilon /4) \frac{1}{c_{*} \pi } \log N + \left\vert  - \frac{1}{c_{*} \pi } \log \epsilon + c_{14} \right\vert }{ \frac{1}{c_{*} \pi } \log N - \left\vert  - \frac{1}{c_{*} \pi } \log \epsilon + c_{14} \right\vert }\\
&\leq & \frac{ \epsilon }{2} + \frac{ \left\vert  - \frac{1}{c_{*} \pi } \log \epsilon + c_{14} \right\vert }{ (1/2) \frac{1}{c_{*} \pi } \log N } \leq \frac{ 3 \epsilon }{4}
\end{eqnarray*}
for $ N \geq \max \{ N_{4}, (4/ \epsilon)(N_{2}+1) \}$ by (\ref{eq:CON1}).
The above two inequalities yield (\ref{eq:CON3}). 
\qed
\vspace{.2cm}

(D){\bf Proof of (\ref{comsi;ex1;eq0001}).}
We show the lower bound of Example \ref{example1}.
In this case, $a=1$, $\beta = 1$, $c_{*} = \frac{\pi}{2a}$ and $\varepsilon = \frac{1}{2}$.
We have $\vert E_{N} \vert = 2c_{11}$ by (\ref{eq:spn4}).
We note
\begin{eqnarray*}
c_{11} = \frac{2^{2+\varepsilon -\beta}\pi^{\varepsilon - \beta}C_{*}}{c_{*}^{2}} = 2^{7/2}\pi^{-5/2}C_{*}.
\end{eqnarray*}
We can choose $C_{*} = 1.225$ by (\ref{eq:VPT}).
So we have
\begin{eqnarray*}
\vert E_{N} \vert \leq 2c_{11} \fallingdotseq 0.633807. 
\end{eqnarray*}
We have
\begin{eqnarray*}
\tilde{\Phi}_{N} \leq \frac{2}{\pi^{2}}\log N + \frac{2}{\pi^{2}}
\end{eqnarray*}
by (\ref{eq:spn5}).
So we can show that
\begin{eqnarray*}
&&\sum_{i=0}^{N-1}p_{i}^{(N)} \leq 1 + \tilde{\Phi}_{N} + \vert E_{N} \vert + \frac{1}{1-\rho_{*}}\\
&&\leq 1 + \frac{2a}{\pi^{2}}\log N + \frac{2}{\pi^{2}} + 0.633807 + \frac{1}{1-\rho_{*}}
\end{eqnarray*}
by (\ref{eq:spn1}), (\ref{eq:spn2}) and (\ref{eq:spn3}).
So we have 
\begin{eqnarray*}
\frac{1}{\sum_{i=0}^{N-1}p_{i}^{(N)}} \geq \frac{1}{1 + \frac{2}{\pi^{2}}\log N + \frac{2}{\pi^{2}} + 0.633807 + \frac{1}{1-\rho_{*}}}
\end{eqnarray*}
by Proposition \ref{bound;thm2}.
It is easily to check $r_{*} \fallingdotseq 0.212207$ (by (\ref{eq:C*r*})) and $\max_{r_*\le|\theta|\le\pi}|\phi(\theta)| \le 0.785802$, then we set $\rho_*= 0.785802$.
Then,
\begin{eqnarray*}
\frac{1}{\sum_{i=0}^{N-1}p_{i}^{(N)}} \geq \frac{1}{\frac{2}{\pi^{2}}\log N + \frac{2}{\pi^{2}} + 6.50503}.
\end{eqnarray*}
So we have (\ref{comsi;ex1;eq0001}).
\qed

\begin{backmatter}
\end{backmatter}

\end{document}